\documentclass[11pt,letterpaper]{amsart}
\usepackage[utf8]{inputenc}
\usepackage{verbatim}
\usepackage{amsmath,amssymb,amsfonts,amsthm}
\usepackage{bigints}
\usepackage{enumitem}
\usepackage{hyperref}
\usepackage{array}

\numberwithin{equation}{section}

\newtheorem{definition}{Definition}[section]

\newtheorem{proposition}[definition]{Proposition}
\newtheorem{theorem}[definition]{Theorem}
\newtheorem{fact}[definition]{Fact}
\newtheorem{lemma}[definition]{Lemma}
\newtheorem{corollary}[definition]{Corollary}
\newtheorem{remark}[definition]{Remark}

\newcommand\R{\mathbb{R}}
\newcommand\N{\mathbb{N}}
\newcommand\Z{\mathbb{Z}}

\newcommand{\norm}[1]{ \left\lVert#1\right\rVert}
\newcommand{\mathbbm}[1]{\text{\usefont{U}{bbm}{m}{n}#1}}

\renewcommand{\Im}{\operatorname{Im}}
\author{Miquel Saucedo}
\address{M.  Saucedo,  Centre de Recerca Matemàtica\\
Campus de Bellaterra, Edifici C
08193 Bellaterra (Barcelona), Spain}
\email{miquelsaucedo98@gmail.com }

\author{Sergey Tikhonov}
\address{S. Tikhonov, 
ICREA, Pg. Lluís Companys 23, 08010 Barcelona, Spain;
Centre de Recerca Matem\`{a}tica\\
Campus de Bellaterra, Edifici C
08193 Bellaterra (Barcelona), Spain;
 and Universitat Autònoma de Barcelona.}
\email{stikhonov@crm.cat}

\title[Extremizers of a Fourier uncertainty principle]{
Extremizers of  a Fourier uncertainty principle related to averaging
}

\begin{document}



\thanks{M. Saucedo is supported by  the Spanish Ministry of Universities through the FPU contract FPU21/04230 and PID2023-150984NB-I00.
 S. Tikhonov is supported by
PID2023-150984NB-I00 and 2021 SGR 00087. This work is supported by
the CERCA Programme of the Generalitat de Catalunya and the Severo Ochoa and Mar\'ia de Maeztu
Program for Centers and Units of Excellence in R\&D (CEX2020-001084-M)}

\begin{abstract}
   We study
   the uncertainty principle 
   $$
\norm{\widehat{\mu}(\xi) |\xi|^\beta}_\infty^{\alpha}
   \left(\int |x|^\alpha d \mu\right)^{\beta}
\geq
C(\alpha,\beta,d){\norm{\mu}_{TV}^{\alpha+\beta}}$$
for finite non-negative measures on $\R^d $. We prove that $C(\alpha,\beta,d)>0$ for all $\alpha,\beta>0$ and that  extremizers exist. Moreover, we obtain an abstract characterization of the extremizers, which allows us 
to describe their asymptotic behavior and, for certain parameter values, to determine them explicitly.
\end{abstract}
\maketitle
\section{Introduction and main results}

In a recent work, Steinerberger \cite{steinerberger} proposed an axiomatic approach to find the \textit{optimal} way of averaging a function over a given scale. He assumed that such an averaging  operator must take the form
$$\operatorname{Av}_f(F):=F*f,$$ where,
for a given $\alpha>0$,
the function $f$ satisfies
$$f\geq 0, \quad \norm{f}_1=1 \mbox{ and }\norm{f(x)|x|^\alpha}_1=1.$$ 
To select the best averaging function, Steinerberger suggested that the average $\operatorname{Av}_f(F)$ should be as smooth as possible. To quantify this smoothness, he considered the quantity 
\begin{equation}
    \label{eq:first}
    \sup_{\norm{F}
_2=1} \norm{\nabla^\beta  \operatorname{Av}_f(F)}_2, \quad \mbox{ for a given } \beta>0.
    \end{equation}
He showed that  this expression  is  minimized when 
 $f$ is an extremizer of the following uncertainty principle:
\begin{equation}
    \label{eq:unc}
    \big\|{\widehat{f}(\xi) |\xi|^\beta}\big\|_\infty^{\alpha} \big\|{f(x)|x|^\alpha}\big\|_1^{\beta} \geq C_{\alpha,\beta,d} \big\|{f}(x)\big\|_1^{\alpha+\beta}, \;\; f\geq 0,\; f\in L^1(\R^d).
    \end{equation}

Motivated by this observation, Steinerberger \cite{steinerberger} studied several properties of \eqref{eq:unc}. 

First, he showed that for
  $\alpha>0$ and $\beta>d/2$, there exists a positive constant $\widetilde{C}_{\alpha,\beta,d}$ such that for all $f\in L^1(\R^d)$ (not necessarily non-negative)
\begin{equation}
    \label{eq:unc2}
    \big\|{\widehat{f}(\xi) |\xi|^\beta}\big\|_\infty^{\alpha} \big\|{f(x)|x|^\alpha}\big\|_1^{\beta} \geq \widetilde{C}_{\alpha,\beta,d} \big\|{f}(x)\big\|_1^{\alpha+\beta}.
    \end{equation}  

Second, as a step towards finding the best averaging function $f$,   Steinerberger
 proved that 
 the characteristic function of the interval is a local extremizer of \eqref{eq:unc} for $d=\beta=1$ and  $\alpha
 \in \{2, 3, 4, 5, 6\}
 $.
He also 
 conjectured that for $(\alpha,\beta,d)=(\alpha,1,1)$ with $\alpha\geq 2$, the characteristic function of the interval is a global extremizer of \eqref{eq:unc}.
 Later, in joint work with Kravitz \cite{steinerberger2}, he proved that $f(x) = 1-|x|$ is a local extremizer for $(\alpha,\beta,d)=(2,2,1)$ in the class of all symmetric $C^3$-functions
compactly supported on $[-1, 1]$ with nonnegative Fourier transform.
Finally, we would like to mention the recent paper
\cite{Richardson}, which studies a discrete version of the averaging problem in the case $\beta=2$.

 The aim of this paper is threefold. First, we establish that the Fourier uncertainty principle
\eqref{eq:unc}
holds for all $\alpha,\beta>0$, and that extremizers exist in every case. Second, we provide an abstract characterization of these extremizers, which allows us to derive the asymptotic behavior of their Fourier transforms. Third, for specific parameter values, we explicitly identify the extremizers by refining the argument for local extremality presented in \cite{steinerberger}. 
In particular, we 
prove that
in the case $(\alpha,1,1)$ with $\alpha\geq 2$
the standard average
$\displaystyle\operatorname{Av}(F)(x):=\int_{x-\frac{1}{2}}^ {x+\frac{1}{2}}F$ is optimal in \eqref{eq:first}, 
thus
confirming Steinerberger’s conjecture.

Before stating our results, we note that for the averaging problem it is  natural to consider  any average
of the form
$$\operatorname{Av}(F)_{\mu}(x)=\int F(y)d \mu(x-y),$$ where $\mu$ is a Borel probability 
 measure. 

\vspace{0.4cm}
{\it{Notation and definitions.}}

The $d$-dimensional Fourier transform of a function $f\in L^1(\R^d)$ is given by
 $$ \widehat{f}(\xi)=\int_{\R^d} {f}(x) e^{-2\pi i \langle x, \xi \rangle} d x.
$$
Let $\mathcal{M}(\R^d)$ denote the space of real-valued, finite Borel measures.
The subspace $\mathcal{M}^+(\R^d)\subset \mathcal{M}(\R^d)$ consists of
 non-negative  
 measures.
 
 For a measure $\mu $ and a measurable function $f$, we use the notation
    $$I_\mu(f):=\int f d \mu.$$

    The variation of a signed measure $\mu$ is
    $$|\mu|=\mu^++ \mu^-,$$ where $\mu^+$ and $\mu^-$ are the positive measures given by the Hahn decomposition theorem. The total variation of $\mu$ is defined by
    $$\norm{\mu}_{TV}= \int d |\mu|.$$

The Fourier transform of a tempered distribution $\omega$ is defined by
$$ \widehat{\omega} (g) = \omega (\widehat{g}) \quad \mbox{ for } g\in S(\R^d). $$ As usual, we may write
$\int \omega g$  instead of $\omega(g)$.
    Observe that if $\norm{\mu}_{TV}<\infty$, then $\widehat{\mu}$ is a bounded, continuous function given by
    $$\widehat{\mu}(\xi)=\int e^{-2 \pi i \langle x ,\xi\rangle} d \mu(x).$$

    We say that $\mu$ is radial if 
    $$\int f d \mu=\int Rf d \mu, \quad\mbox{ for } f\in C(\R^d) \mbox{ and } R\in O({d}),$$ where
    $$Rf(x)=f(Rx).$$ Observe if $\mu$ is real and radial, then $\widehat{\mu}$ is also real and radial.

We write $\lambda \mathbb{S}^{d-1}:=\{|x|=\lambda\}$ for the sphere of radius $\lambda$ and denote by $\delta_{\lambda \mathbb{S}^{d-1}}$ the uniform uniform measure on this sphere, normalized so that $\int_{\R^d} \delta_{\lambda \mathbb{S}^{d-1}}=\lambda^{d-1} \omega_{d-1},$ where $\omega_{d-1}$ is the surface measure of $\mathbb{S}^{d-1}$.
  \begin{definition}Let $\alpha, \beta>0$.
Let $\mu \in \mathcal{M}^+(\mathbb{R}^d)$ and suppose that 
\[
\|\widehat{\mu}(\xi) |\xi|^\beta\|_\infty < \infty.
\]
We define
\[
C_{\alpha,\beta,d}(\mu) := 
\frac{\|\widehat{\mu}(\xi) |\xi|^\beta\|_\infty^{\alpha} I_\mu(|x|^\alpha)^{\beta}}{\|\mu\|_{\mathrm{TV}}^{\alpha + \beta}}.
\]
If the above norm is infinite, we set \( C_{\alpha,\beta,d}(\mu) := \infty \).

We then define the optimal constant
\[
C^*_{\alpha,\beta,d} := \inf_{0 \neq \mu \in \mathcal{M}^+(\mathbb{R}^d)} C_{\alpha,\beta,d}(\mu).
\]

When \( d\mu = f\,dx \), we will write \( C_{\alpha,\beta,d}(f) \) instead of \( C_{\alpha,\beta,d}(\mu) \).
\end{definition}

With this notation, the extension of inequality~\eqref{eq:unc} to measures becomes
\begin{equation} \label{uncer:measures}
\|\widehat{\mu}(\xi)|\xi|^\beta\|_\infty^\alpha\, I_\mu(|x|^\alpha)^\beta 
\geq C^*_{\alpha,\beta,d}\|\mu\|_{\mathrm{TV}}^{\alpha + \beta}.
\end{equation}

\vspace{0.4cm}
{\it{Main results.}}

Before formulating the main results of the paper, we 
 demonstrate  that the uncertainty principle  for  finite non-negative measures always holds  and admits an extremizer.

\begin{theorem} 
\label{theorem:notmain} Let $\alpha, \beta>0$. Then there exists a radial non-zero measure $\mu^*\in \mathcal{M}^+(\R^d)$ such that
 $C_{\alpha,\beta,d}(\mu^*)=C^*_{\alpha,\beta,d}>0$.
\end{theorem} 
\begin{remark}
\label{remark:first}
In \cite{steinerberger} it was shown that  $\inf_{f \neq 0}C_{\alpha,\beta,d}(f)>0$ for $\beta>\frac d2$. This is optimal, since if $\beta\leq \frac d2 $, then $\inf_{f \neq 0}C_{\alpha,\beta,d}(f)=0$.
\end{remark}

Our first main result identifies the extremizers of 
\eqref{uncer:measures}
 for specific  triples  $(\alpha,\beta,d)$.

\begin{theorem} 
\label{theorem:main}
    \begin{enumerate}
  The following hold:     
        \item
        Setting
$f_1(x)=\mathbbm{1}_{[-\frac12,\frac12]}(x),$ 
        for any $\alpha\geq 2$,$$C^*_{\alpha,1,1}=C_{\alpha,1,1}(f_1)=\frac{1}{(2\pi)^\alpha(\alpha+1)};$$  
        \item  Setting
        $f_2(x)=(1-|x|)\mathbbm{1}_{[-1,1]}(x)$,
        $$C^*_{1,2,1}=C_{1,2,1}(f_2)=\frac{1}{9 \pi^2};$$ 
         \item Setting
        $f_3(x)=|x|^{-1}\mathbbm{1}_{[0,1]}(|x|)$,
$$C^*_{1,2,3}=C_{1,2,3}(f_3)=\frac{4}{9 \pi ^2};$$
        \item Setting $\mu_4= \delta_{\mathbb{S}^2}$, 
        $$C^*_{\alpha,\beta,3}= C_{\alpha,\beta,3}(\mu_4)=\left(\frac{\sin{2 \pi \lambda_\beta}}{2\pi \lambda_\beta^{1- \beta}}\right)^\alpha,$$ provided that $$(\alpha,\beta) \in \Big\{(\alpha,1): \alpha \geq  ({\frac{\pi}{2}-1})^{-1}\Big\} \cup \bigg\{(\alpha,\beta): \alpha\geq 2, 0<\beta <1\bigg\},$$ where $\lambda_\beta$ is defined implicitly by
    \begin{equation*}
    \label{eq:deflb2}
     2 \pi \lambda_\beta =(1-\beta)\tan{2 \pi \lambda_\beta }, \; \lambda_\beta \in (0,\frac14),\; \beta<1, 
    \end{equation*} and $\lambda_1=\frac14$.
    \end{enumerate}
\end{theorem}

Our second main result describes 
the general structure of extremizers.
It turns out that there are two possibilities, depending on the value of the parameter $\beta$.
\begin{theorem}\label{extremizer}
    Let $\alpha,\beta>0$. 
    Assume that the radial measure $\mu^* \in \mathcal{M}^+(\R^d)$ satisfies $C_{\alpha,\beta,d}(\mu^*)=C^*_{\alpha,\beta,d}.$ We have

    \begin{itemize}
        \item if $\beta>\frac{d-1}{2}$,
        then
        for any $\Delta>0$ there exists $|\xi|>\Delta$ such that 
        \begin{equation}
    \label{mu-norm at infinity}
    \norm{\widehat{\mu^*}(\xi)|\xi|^\beta}_\infty=|\widehat{\mu^*}(\xi)||\xi|^\beta;
        \end{equation}

\item if $\beta\leq \frac{d-1}{2}$, then
$\mu^*$ is compactly supported and 
either  the above holds 
 or      there exist $N,c_k,\lambda_k\geq 0$ such that \begin{equation}
    \label{mu* is finite}
    \mu^*=\sum_{k=0}^N c_k\delta_{\lambda_k \mathbb{S}^{d-1}}.
    \end{equation}
\end{itemize}

    \end{theorem}

Our next result determines the decay of the Fourier transform of the extremizers of \eqref{uncer:measures}, answering a question posed in \cite{steinerberger}.

\begin{corollary}
 \label{th:asymptotic}
   Let $\alpha,\beta>0$. Assume that the radial measure $\mu^* \in \mathcal{M}^+(\R^d)$ satisfies $C_{\alpha,\beta,d}(\mu^*)=C^*_{\alpha,\beta,d}.$ Then
    $$0<\limsup_{|\xi|\to \infty} |\widehat{\mu^*}(\xi) |\xi|^\gamma<\infty,$$ where $\gamma= \beta$ if        $\beta\geq \frac{d -1}{2}$ and $\gamma \in \{ \beta,\frac{d -1}{2}\}$ otherwise.
\end{corollary}

\begin{remark}
    \label{remark:big} The following observations complement Theorems \ref{theorem:main},\ref{extremizer} and Corollary \ref{th:asymptotic}:
    \begin{enumerate}[label=(\roman*)]
               \item In Theorem \ref{theorem:main}, $f_1$ is the unique extremizer up to dilation and normalization.  In cases (2)–(4), the presented extremizers  are unique among radial functions, but non-radial extremizers also exist. For example, $\tilde{f_2}(x)=f_2(x) + \varepsilon \mathbbm{1}_{[-1,1]}(x) \sin (2 \pi x)$ remains an extremizer for $ |\varepsilon|$ small enough.
        \item The measure $\mu_4$ appears to be an extremizer for other values of $\alpha,\beta$ as well. Verifying this could reduce to proving a variant of Lemma \ref{lemma:b1}
        \item   Regarding the case $\beta<\frac{d-1}{2}$ in
        Theorem \ref{extremizer} and Corollary \ref{th:asymptotic}, 
        from the extremality of $\mu_4$ we see that the optimal $\mu^*$ can have the representation \eqref{mu* is finite}. Moreover, 
        $\widehat{\mu_4}(\xi)$  decays like $|\xi|^{-\frac{d-1}{2}}$, rather than the expected $|\xi|^{-\beta}$. We do not know if this is always the case for $\beta<\frac{d-1}{2}$.
    \end{enumerate}
\end{remark}

\vspace{0.4cm}
{\it{Structure of the paper.}}
The rest of the paper is organized as follows.

In Section \ref{sec2}, we prove that the uncertainty principle
\eqref{uncer:measures}
holds for all finite non-negative measures and for any $\alpha,\beta>0$ and  establish the existence of extremizers.We also demonstrate Remark \ref{remark:first}.

Section \ref{sec3} is devoted to the proof of Theorem \ref{theorem:main} and Remark \ref{remark:big}(i). Our main tool is Lemma \ref{lemma:dual}, where we introduce a general framework for verifying extremality. This result can be interpreted as an instance of weak duality in linear programming (or convex optimization).

In Section \ref{sect:duality}, we  develop the duality technique to prove Theorem \ref{extremizer} and Corollary \ref{th:asymptotic}.

Section \ref{sec5} contains additional applications of the methods developed in the earlier sections. In particular, we address a question posed by Kravitz and Steinerberger in \cite{steinerberger2}, confirming that the function
$
f(x)=(1-|x|)\mathbbm{1}_{[-1,1]}(x)$ is a global extremizer of the inequality
$$
    \big\|{\widehat{f}(\xi) |\xi|^2}\big\|_\infty^{2} \big\|{f(x)|x|^2}\big\|_1^{2} \geq {C}^*_{2,2,1} \big\|{f}(x)\big\|_1^{4}
  $$ 
  among all
positive-definite functions.

In the Appendix, we prove some auxiliary inequalities involving trigonometric functions, which we need for the proof of Theorem \ref{theorem:main}.
\vspace{0.3cm}

\section{Proof of Theorem \ref{theorem:notmain} }\label{sec2}
    \begin{proof}[Proof of Theorem \ref{theorem:notmain}]
    To begin with, observe that by a scaling argument it suffices to consider $\mu\in \mathcal{M}^+(\R^d)$  normalized so that $\norm{\mu}_{TV}=1$ and $\norm{\widehat{\mu}(\xi)|\xi|^\beta}_\infty=1$. Let then $(\mu_n)_{n=1}^\infty $ be a sequence of normalized probability measures such that $$C_{\alpha,\beta,d}(\mu_n)\searrow C^*_{\alpha,\beta,d}\geq 0.$$

   Observe that the  sequence $(\mu_n)_{n=1}^\infty$ is tight. Indeed, let $\varepsilon>0$, then for every $n\geq 1$,
    $$\int_{|x|>\varepsilon^{-1}} d\mu_n \leq \varepsilon^\alpha I_{\mu_n}(|x|^\alpha) = \varepsilon^\alpha C^{\frac1 \beta}_{\alpha,\beta,d}(\mu_n)\leq \varepsilon^\alpha C^{\frac1 \beta}_{\alpha,\beta,d}(\mu_1).  $$

    Hence, by Prokhorov's Theorem \cite[Theorem 5.1]{prokh}, there exist a probability measure $\mu^*$ and a subsequence $(\mu_{n_k})_{k=1}^\infty$ such that for all bounded continuous functions $f$

    $$\int f d\mu_{n_k} \to \int f d\mu^* .$$
    In particular, $\widehat{\mu^*}(\xi)=\lim_k \widehat{\mu_{n_k}}(\xi)$, whence we deduce that $\norm{\widehat{\mu^*}(\xi)|\xi|^\beta}_ \infty \leq 1. $ Moreover, for any compactly supported continuous function $\varphi\leq 1$,
    $$\int \varphi(x) |x|^\alpha d\mu^*(x) = \lim_{k} \int \varphi(x) |x|^\alpha d\mu_{n_k} \leq \lim_{k} I_{\mu_{n_k}}(|x|^\alpha) = C^{*{\frac{1}{\beta}}}_{\alpha, \beta,d}, $$
    which implies that $$C_{\alpha,\beta,d}(\mu^*) =\norm{\widehat{\mu^*}(\xi)|\xi|^\beta}^\alpha_ \infty I^\beta_{\mu^*}(|x|^\alpha) \leq C^{*}_{\alpha, \beta,d}$$ and $C_{\alpha,\beta,d}(\mu^*)= C^{*}_{\alpha, \beta,d}.$

    Second, to show that $C^*_{\alpha,\beta,d}>0$, observe that if $C_{\alpha,\beta,d}(\mu^*)=0$, then either $\widehat{\mu^*}(\xi)=0$ for $\xi\neq 0$ or $I_{\mu^*}(|x|^\alpha)=0$. The former is impossible because $\mu^*\neq 0$ and $\widehat{\mu^*}$  is continuous. If the latter holds, then $\mu^*=\delta_0$, which contradicts the fact that $\norm{\widehat{\mu^*}(\xi)|\xi|^\beta}\leq 1$. In conclusion, $C_{\alpha,\beta,d}(\mu^*)=C^*_{\alpha,\beta,d}>0.$

    Finally, to see that $\mu^*$ can be taken to be radial, we define the radial measure $\mu^*_{rad}$ by 
    $$\int f d \mu^*_{rad}= \int f_{rad} \;d \mu^* , \mbox{ for } f \in C(\R^d), $$ where
    $$f_{rad}(x)=  \frac{1}{\omega_{d-1}}\int_{\mathbb{S}^{d-1}} f(\nu |x|)d \nu.$$
    Note that $I_{\mu^*}(|x|^\alpha)= I_{\mu^*_{rad}}(|x|^\alpha)$, $\norm{\mu^*}_{TV}=\norm{\mu^*_{rad}}_{TV}$ and
    $$\widehat{\mu^*_{rad}}(\xi) = \frac{1}{\omega_{d-1}}\int_{\mathbb{S}^{d-1}} \widehat{\mu^*}(\nu |\xi|)d \nu.$$ Therefore, 
    $$C_{\alpha,\beta,d}(\mu^*_{rad})\leq C_{\alpha,\beta,d}(\mu^*), $$ whence the result follows.

    \end{proof}

\begin{proof}[Proof of Remark \ref{remark:first}]
Let $\widehat{\Phi}$ be a $C^\infty(\R^d)$ function whose support is contained in the unit square. We know that $\Phi$ belongs to the Schwartz class of rapidly decreasing functions.

For $M \in \mathbb{N}$ and $(c_j)_{j\in \mathbb{N}^d}$ a compactly supported sequence to be fixed, consider the following function: $$f(x)= \Phi(x) \sum_{j \in \mathbb{N}^d}  \varepsilon_{j} c_{j} e^{-2 \pi i \langle x,j \rangle},$$ where each $\varepsilon_j$ is either $+1$ or $-1$.

First,
by using Fubini and Khinchin's inequality, averaging over all possible values for $\varepsilon$, we estimate $\mathbbm{E}\norm{f}_1$ as follows:
\begin{align*}
\mathbbm{E}\left[\int_{\R^d} \Big|\Phi(x) \sum_{j \in \mathbb{N}^d}  \varepsilon_{j} c_{j} e^{-2 \pi i \langle x,j \rangle} \Big| dx \right]&=\int_{\R^d} |\Phi(x) | \mathbbm{E}\left[\Big|\sum_{j\in \mathbb{N}^d}  \varepsilon_{j} c_{j} e^{-2 \pi i \langle x, j\rangle } \Big|\right]dx \\
&\approx \int_{\R^d} \left| \Phi(x) \Big(\sum_{j\in \mathbb{N}^d}   |c_{j}|^2\Big)^{\frac{1}{2}} \right| dx\\
&\approx  \Big(\sum_{j\in \mathbb{N}^d}   |c_{j}|^2\Big)^{\frac{1}{2}}.
\end{align*}
Similarly,
$$\mathbbm{E}\left[\int_{\R^d} \Big||x|^{\alpha}f (x)\Big|dx\right] \approx \left(\sum_{j\in \mathbb{N}^d}   |c_{j}|^2\right)^{\frac{1}{2}}.$$

Second, a straightforward computation shows that $$\widehat{f}(\xi)= \sum_{j\in \mathbb{N}^d}  \varepsilon _{j} c_{j} \widehat{\Phi}(\xi +j).$$ Thus, since the supports of the translates of $\widehat{\Phi}$ are pairwise disjoint,
$$\norm{\widehat{f} (\xi)|\xi|^{\beta}}_{\infty}\approx \left(\sup_{j\in \mathbb{N}^d} \,(|j|+1)^{\beta }|c_{j}|^{ }\right).$$ 

Finally, if the inequality
$$\norm{f}_{1}\lesssim \norm{f (x)|x|^{\alpha}}_1^{\frac{\beta}{\alpha + \beta}}\norm{\widehat{f}(\xi) |\xi|^{\beta}}_{\infty }^{\frac{\alpha}{\alpha + \beta}}$$
holds,
taking expected values and noting that for any  random variable $X\geq 0$ and $0<p\leq 1$ Jensen's inequality gives $$\mathbbm{E}[X^p] \leq \mathbbm{E}[X]^p,$$ we obtain
$$\left(\sum_{j\in \mathbb{N}^d}   |c_{j}|^2\right)^{\frac{1}{2}}\lesssim \left(\sum_{j\in \mathbb{N}^d}   |c_{j}|^2\right)^{\frac{\beta}{\alpha + \beta}\frac{1}{2}}
\left(\sup_{j\in \mathbb{N}^d} (|j|+1)^{\beta }|c_{j}|^{ }\right)^{\frac{\alpha}{\alpha + \beta}}$$ for any sequence $(c_j)$.
The latter is only true if
$\beta>\frac d2.$
\end{proof}

\vspace{0.3cm}
\section{Proofs of Theorem \ref{theorem:main}
 and Remark \ref{remark:big}(i)}\label{sec3}
In the case $(\alpha,\beta,d)=(2,1,1)$, Steinerberger's proof of the local extremality of $f_1$ is based on the following remarkable fact:
\begin{fact}\label{fact}
    For any even $f:[-\frac12, \frac12] \to \R$,
    $$\int_{-\frac12}^\frac12 f(x) \left(1-4x^2 \right) dx \leq \frac{2 \pi}{3} \norm{\widehat{f}(\xi) |\xi|^\beta}_\infty,$$
    with equality for $f=f_1:=\mathbbm{1}_{[-\frac12, \frac12]}$.
\end{fact}
From this result, it is straightforward to deduce that 
$f_1$ is a global extremizer of \eqref{eq:unc} with $C^*_{2,1,1}=\frac{1}{12 \pi^2}$
among the set of nonnegative functions supported on $[-\frac12, \frac12]$ and which satisfy the normalization conditions $$\norm{f}_1= \norm{f_1}_1=1 \mbox{ and } \norm{\widehat{f}(\xi) |\xi|^\beta }_\infty= \norm{\widehat{f_1}(\xi) |\xi|^\beta}_\infty=\frac{1}{\pi}.$$ Indeed, for any such $f$ we have
$$0\leq \int_{-\frac12}^\frac12 f(x) \left(1-4x^2 \right) dx=\norm{f}_1-4 \norm{f(x) |x|^2}_1 \leq \frac{2 \pi}{3} \norm{\widehat{f}(\xi) |\xi|^\beta}_\infty, $$ and using the normalization conditions and rearranging we arrive at

$$4\norm{f_1(x)|x|^2}_1= \norm{f_1}_1-  \frac{2 \pi}{3} \norm{\widehat{f_1}(\xi) |\xi|^\beta}_\infty  \leq  4 \norm{f(x) |x|^2}_1.$$

\vspace{3mm}
The main tool to prove Theorem \ref{theorem:main} is Lemma \ref{lemma:dual},
which can be seen as 
 a further elaboration on the  preceding idea.
 We also remark that both Fact \ref{fact} and Lemma \ref{lemma:dual} can  be understood as a type of linear programming weak duality. This connection with linear programming will be explored further in Section \ref{sect:duality}.

\begin{lemma}
\label{lemma:dual}
    Let $\mu^*\in \mathcal{M}^+(\R^d)$ be radial and  satisfy $$I_{\mu^*} (|x|^\alpha)+ \norm{\widehat{\mu^*}(\xi)|\xi|^\beta} <\infty.$$  Define $$\Lambda=\Big\{|\xi|: \widehat{\mu^*}(\xi)|\xi|^\beta =\pm\norm{\widehat{\mu^*}(\xi)|\xi|^\beta} \Big\}\subset \R_+.$$ 
    Let \begin{equation}
        \label{eq:tekni}\psi=\sum_{\lambda\in \Lambda} c_\lambda \widehat{\delta_{\lambda \mathbb{S}^{d-1}}}\quad \mbox{ with } \quad\sum_{\lambda\in \Lambda} |c_\lambda| \lambda^{d-1}<\infty.
    \end{equation}
    Assume further that
    \begin{enumerate}[label=(\roman*)]
    \item for some $C,D\in \R$ with $C>0$ the function $$H(x):=\psi(x)+C|x|^\alpha-D$$ is non-negative;
        \item $H d \mu^*\equiv 0$;
        \item for every $|\xi|=\lambda\in \Lambda$, we have  $c_\lambda \widehat{\mu^*}(\xi)\geq 0$.
    \end{enumerate} Then $$C_{\alpha,\beta,d}(\mu^*)=C^*_{\alpha,\beta,d}.$$ Moreover, any  $\mu\in \mathcal{M}^+(\R^d)$  such that $\norm{\mu}_{TV}=\norm{\mu^*}_{TV}$ and $\norm{\widehat{\mu}(\xi)|\xi|^\beta}_\infty=\norm{\widehat{\mu^*}(\xi)|\xi|^\beta}$ satisfies $C_{\alpha,\beta,d}(\mu)=C^*_{\alpha,\beta,d}$  if and only if 
    \begin{enumerate}[label=(\alph*)]
        \item $Hd\mu\equiv 0$;
        \item for every $|\xi|=\lambda\in \Lambda$, we have $c_\lambda |\xi|^\beta \widehat{\mu}(\xi) = |c_\lambda| \norm{\widehat{\mu}(\xi)|\xi|^\beta}_\infty.$
    \end{enumerate}
\end{lemma}
\begin{proof}

    Let $\mu\in \mathcal{M}^+(\R^d)$ satisfy $$\norm{\mu^*}_{TV}=\norm{\mu}_{TV}, \norm{\widehat{\mu^*}(\xi) |\xi|^\beta}_\infty=\norm{\widehat{\mu}(\xi) |\xi|^\beta}_\infty \mbox{ and } I_{\mu}(|x|^\alpha)<\infty.$$ The result follows if we prove that $I_{\mu^*}(|x|^\alpha)\leq I_{\mu}(|x|^\alpha)$. To see this, first note that by Plancherel and  \eqref{eq:tekni}
\begin{equation}
\label{eq:gg}
\begin{aligned}
    \int_{\R^d}  \psi d\mu &=  \sum_{\lambda\in \Lambda} c_\lambda\int_{\R^d}    \widehat{\delta_{\lambda \mathbb{S}^{d-1}}}d\mu\\
    &=\sum_{\lambda\in \Lambda} c_\lambda \int_{\lambda \mathbb{S}^{d-1}} \widehat{\mu}   \leq \omega_{d-1}\norm{\widehat{\mu}(\xi)|\xi|^\beta}_\infty \sum_{\lambda\in \Lambda} \lambda ^{d-1-\beta}|c_\lambda|.
    \end{aligned}
\end{equation}For $\mu^*$, because of property (iii) and the definition of $\Lambda$, we actually have
\begin{equation}
\label{eq:gf}
    \int_{\R^d}  \psi d \mu^* =\omega_{d-1}\sum_{\lambda\in \Lambda} \lambda ^{d-1}|c_\lambda| |\widehat{\mu^*}(\lambda)|=\omega_{d-1}\norm{\widehat{\mu^*}(\xi)|\xi|^\beta} \sum_{\lambda\in \Lambda} \lambda ^{d-1-\beta}|c_\lambda| .
\end{equation}

    Thus, letting 
   $A=\omega_{d-1} \sum_{\lambda\in \Lambda} \lambda ^{d-1-\beta}|c_\lambda|$, we have, in light  of property (ii) and \eqref{eq:gf}, \begin{eqnarray*}
    0 &=& \int_{\R^d}  Hd \mu^*= \int_{\R^d}  \psi d \mu^* +CI_{\mu^*}(|x|^\alpha) - D\norm{\mu^*}_{TV}\\
    &=& A\norm{\widehat{\mu^*}(\xi) |\xi|^\beta}_\infty  +CI_{\mu^*}(|x|^\alpha) - D\norm{\mu^*}_{TV}.
    \end{eqnarray*}
  For $\mu$, using now \eqref{eq:gg}, we see that
    \begin{equation}
    \begin{aligned}
    \label{eq:ggg}
    0 &\leq \int_{\R^d}{H} d \mu  = \int_{\R^d}  \psi d \mu +CI_{\mu}(|x|^\alpha) - D\norm{\mu}_{TV}\\
&\leq A\norm{\widehat{\mu}(\xi) |\xi|^\beta}_\infty  +CI_{\mu}(|x|^\alpha) - D\norm{\mu}_{TV}\\
    &= A\norm{\widehat{\mu^*}(\xi) |\xi|^\beta}_\infty  +CI_{\mu}(|x|^\alpha) - D\norm{\mu^*}_{TV},
    \end{aligned}
    \end{equation}
   whence we conclude that $I_{\mu^*}(|x|^\alpha)\leq I_{\mu}(|x|^\alpha).$ 
   
   Moreover, any non-negative measure $\mu$ such that $\norm{\mu}_{TV}=\norm{\mu^*}_{TV}$ and $\norm{\widehat{\mu}(\xi)|\xi|^\beta}_\infty=\norm{\widehat{\mu^*}(\xi)|\xi|^\beta}_\infty$  is an extremizer if and only if all inequalities in \eqref{eq:gg} and \eqref{eq:ggg} are equalities, that is, if and only if $\displaystyle \int H d \mu=0$ and $\displaystyle c_\lambda \lambda^\beta \int_{\lambda \mathbb{S}^{d-1}} \widehat{\mu}= |c_\lambda | \omega_{d-1} \lambda^{d-1} \norm{\widehat{\mu}(\xi)|\xi|^\beta}_\infty $ hold, which are equivalent to (a) and (b).
\end{proof}
We now proceed to prove Theorem \ref{theorem:main} and Remark \ref{remark:big} (i). 

\begin{proof} [Proof of item (1) in Theorem \ref{theorem:main}.]
 Recall that $\alpha \geq 2$ and $\beta=1$.

Set $$f_1=\mathbbm{1}_{[-\frac{1}{2},+\frac{1}{2}]}\qquad \mbox{and
}\qquad H(x)=\psi_1(x)+|2x|^\alpha -1,$$ where   $$\psi_1(x)= 2\sum_{k=0}^\infty c_k^{(\alpha)} \cos{2 \pi  (k+\frac{1}{2})x},\quad 
c_k^{(\alpha)}=\int_{0}^{\frac{1}{2}} (1-|2x|^\alpha) \cos{2 \pi  (k+\frac{1}{2})x} dx.$$ Assume momentarily that
\begin{equation}
\label{eq:coeffs}
     1\gtrsim c^{(\alpha)}_k (-1)^{k} (k+\frac{1}{2})^3 >0.
\end{equation}
We now verify the hypotheses of Lemma \ref{lemma:dual}.

First, since $\widehat{f_1}(\xi) = \frac{\sin(\pi \xi)}{\pi \xi}$ and $\norm{\widehat{f}_1(\xi)\xi}_\infty=\frac{1}{\pi}$, we have that $\Lambda=\mathbb{N}+\frac12$ and $\widehat{f_1}(k+\frac{1}{2})= \frac{1}{\pi}\frac{(-1)^{k}}{k+\frac{1}{2}}.$ Second, for $d=1$, $\widehat{\delta_{\lambda \mathbb{S}^{d-1}}}(x)=2 \cos{2 \pi \lambda x}$. Thus, we indeed have
$$\psi_1=\sum_{\lambda\in \Lambda} c_\lambda \widehat{\delta_{\lambda \mathbb{S}^{d-1}}},$$ with $c_\lambda= c_k^{(\alpha)}$ for $\lambda=k+\frac12$. Moreover, in light of \eqref{eq:coeffs},  $\sum_{\lambda\in \Lambda} |c_\lambda|<\infty$ and $\widehat{f_1}(\lambda) c_\lambda \geq 0$ for $\lambda \in \Lambda$, that is, both \eqref{eq:tekni} and property (iii) hold.

Third, we  prove items (i) and (ii). Noting that $\psi_1(x)=1-|2x|^\alpha$ for $|x|\leq \frac{1}{2}$,  we immediately see that $f_1H\equiv 0$. Then, by using $\psi_1(x+1)=-\psi_1(x)$, we obtain
  $$H(x)=|2x|^{\alpha} + |2(x-1)|^{\alpha} -2> 0, \quad x\in( 1/2,  3/2].$$
  For $x>\frac32$, we deduce that $H(x)>0$ because $|2x|^\alpha-1$ is strictly increasing and $\psi_1$ is $2$-periodic. Finally, for $x<0$ the result follows by evenness. 
Thus, from Lemma \ref{lemma:dual} we deduce that $f_1$ is an extremizer. Here we note that a computation shows that in this case $C^*_{\alpha,1,1}=\frac{1}{(2\pi)^\alpha (\alpha+1)}$.

Besides, up to scaling and dilation, $f_1$ is the unique extremizer. Indeed, by Lemma \ref{lemma:dual}, any other extremizer $\mu$ normalized such that $\norm{\mu}_{TV}=1$ and $\norm{\widehat{\mu}(\xi) |\xi|}_\infty=\frac1\pi$ must be supported on $[-\frac12,\frac12]$ and satisfy $$\widehat{\mu}(-k-\frac12)=\widehat{\mu}(k+\frac12)=\frac{(-1)^k}{\pi(k+\frac12)}=\widehat{f_1}(k+\frac12) = \widehat{f_1}(-k-\frac12) \mbox{ for } k\in \mathbb{N}.$$ 
This implies that $\widehat{f_1}=\widehat{\mu}$ by observing that the entire function $$g(z):=\frac{\widehat{f_1}(z)-\widehat{\mu}(z)}{\cos{ \pi z}}$$ is bounded on the complex plane and satisfies $\liminf_{x\to \infty} |g(x)|= 0$. 

All that remains is to obtain \eqref{eq:coeffs}. (We note that this fact was also obtained in  \cite{cho2021zeroscertainfouriertransformsimprovements}, here we give a simpler proof.)
For $\alpha=2$, the result holds because $$c^{(2)}_k=\frac{16}{\pi^3} \frac{(-1)^{k}}{(2k+1)^3}.$$
For $\alpha>2$, by repeated integration by parts and setting $g(x)=1-|2x|^\alpha$,
\begin{align*}
      c_k^{(\alpha)}&=\frac{1}{2\pi (k+\frac{1}{2})} \int_0 ^{\frac{1}{2}} -g'(x) \sin\left(2 \pi x (k+\frac{1}{2})\right)dx\\
    &= \frac{1}{\left(2\pi (k+\frac{1}{2})\right)^2} \int_0 ^{\frac{1}{2}} -g''(x) \cos\left(2 \pi x (k+\frac{1}{2})\right)dx\\
   &=\frac{1}{\left(2\pi (k+\frac{1}{2})\right)^3} \int_0 ^{\frac{1}{2}} -g'''(x) \left(\sin\left( \pi  (k+\frac{1}{2})\right)-\sin\left(2 \pi x (k+\frac{1}{2})\right)\right)dx,
\end{align*}
    and \eqref{eq:coeffs} follows by noting that $\sin\left( \pi  (k+\frac{1}{2})\right)=(-1)^{k}$ and, as $\alpha>2$, $g'''$ is integrable and $-g'''(x)>0$.

\end{proof}

\begin{proof}[Proof of item (2) in Theorem \ref{theorem:main}]
Recall that $\alpha=1$ and $\beta=2$.

Set $$f_2(x)=(1-|x|)\mathbbm{1}_{[-1,1]}(x)\qquad \mbox{and
}\qquad H(x)=\psi_2(x) + |x| -\frac 12,$$ where   $$\psi_2(x)= \frac{4}{\pi^2} \sum_{k =0}^\infty\frac{\cos{2\pi(k+\frac12) x}}{(2k+1)^2}.$$

We now check that the assumptions of Lemma \ref{lemma:dual} are satisfied.

 First, since  $\widehat{f_2} (\xi)= \left(\frac{\sin(\pi \xi)}{\pi \xi}\right)^2$ and $\norm{\widehat{f}(\xi)|\xi|^2}_\infty=\frac{1}{\pi^2}$, we have that $\Lambda=\N+\frac12$. 
 Thus, we indeed have
$$\psi_2=\sum_{\lambda\in \Lambda} c_\lambda \widehat{\delta_{\lambda \mathbb{S}^{d-1}}},$$ with $c_\lambda= \frac{2}{\pi^2} \frac{1}{(2k+1)^2} $ for $\lambda=k+\frac12$, which means that \eqref{eq:tekni} holds. Second, from the positivity of $\widehat{f}$ and $c_\lambda$ we see that property (iii) of Lemma \ref{lemma:dual} is also true.

Third, to prove items (i) and (ii) in Lemma \ref{lemma:dual}, in light of 
$$\int_{-1}^1 (\frac12-|x|) \cos({ \pi n x}) d x=2\frac{1-(-1)^n}{(\pi n)^2},\quad n\geq 1,$$ and $\displaystyle\int_{-1}^1 \big(\frac12-|x|\big)  d x=0,$ we establish 
$$\frac12- |x|=\sum_{n=1}^\infty 2\frac{1-(-1)^n}{(\pi n)^2} \cos{n\pi x}=\psi_2(x),\quad -1\leq x \leq 1.$$
Therefore, $H\equiv 0$ on $[-1,1]$, which means that $f_2H\equiv 0$. Besides, we see that $H\geq 0$ because, using that $\psi_2$ is 2-periodic, we have $$H(x)=|x| - |x-2k|>0,\; 2k-1< x <2k+1, \; k \in \Z\setminus \{0\}.$$ 

Thus, Lemma \ref{lemma:dual} implies that $f_2$ is an extremizer.
In this case, direct computation yields $C_{1,2,1}=\frac{1}{9 \pi ^2}$.

To see that $f_2$ is the unique {\it even} extremizer, observe that Lemma \ref{lemma:dual} implies that every extremizer $\mu\in \mathcal{M}^+(\R)$ normalized so that $\norm{\mu}_{TV}=\norm{f_2}_1$ and $\norm{\widehat{\mu}(\xi)|\xi|^2}_\infty=\norm{\widehat{f_2}(\xi)|\xi|^2}_\infty=\frac{1}{\pi^2}$ must be supported on $[-1,1]$ and satisfy, by (b) in Lemma \ref{lemma:dual},

\begin{equation}
\label{vspm}
    \widehat{\mu}(k+\frac12)=\widehat{f_2}(k+\frac 12).
\end{equation}
Next, since $\mu$ is even, $\widehat{\mu}$ is real-valued, so using the fact that the smooth function $\xi^2 \widehat{\mu}(\xi)$ has a maximum at $\xi=k+\frac12$ for $k\in \Z$, by computing the derivative we can deduce that
$$\widehat{\mu}'(k+\frac12)=-\frac{2}{k+\frac12} \widehat{\mu}(k+\frac12)=\widehat{f_2}'(k+\frac12), \; k\in \Z,$$ which, together with \eqref{vspm}, is enough to conclude that $\widehat{f_2}=\widehat{\mu}$. This can be seen by observing that the entire function 
$$g(z)=\frac{\widehat{f_2}(z)-\widehat{\mu}(z)}{\cos^2{ \pi z}}$$ is bounded and satisfies $\liminf_{x\to \infty} |g(x)|=0$.

Finally, we show that $f_2$ is not the unique extremizer. For this purpose, set $$\tilde{f_2}:= f_2+ \varepsilon \varphi,$$ where $\varphi(x)= \mathbbm{1}_{[-1,1]}\sin{ 2 \pi x},$ and observe that for $|\varepsilon|>0$ small enough, $\tilde{f_2}\geq 0$. Thus, since $\varphi$ is odd, we have
$$\norm{f_2}_1=\norm{\tilde{f_2}}_1 \mbox{ and } \norm{f_2(x) |x|}_1=\norm{\tilde{f_2}(x) |x|}_1,$$ so to show that $\tilde{f_2}$ is an extremizer it suffices to prove that \begin{equation}
\label{eq:goal}
    \norm{\widehat{f_2}(\xi)|\xi|^2}_\infty =\norm{\widehat{\tilde{f_2}}(\xi)|\xi|^2}_\infty.
\end{equation}

To see this, note that since $\varphi$ is odd, $\widehat{\varphi}$ is purely imaginary, which means that $$|\widehat{\tilde{f_2}}|^2= |\widehat{f_2}|^2+ \varepsilon^2|\widehat{\varphi}|^2.$$ Next, observe that the function $$g(\xi):=\left |\xi^2 \widehat{\varphi}(\xi)\right|^2=\left(\frac{\xi^2}{\pi(1-\xi^2)}\right)^2\sin^2{2 \pi \xi}$$ satisfies the following properties:

\begin{enumerate}
    \item $
        g(k+\frac12)=0, \; k \in \Z;$
   
    \item $        g'(k+\frac12)=0, \; k \in \Z;$
     \item $\norm{g}_\infty+\norm{g''}_\infty<\infty$.
    
\end{enumerate}
In light of property (1), to obtain \eqref{eq:goal}, it suffices to show that for $|\varepsilon|$ small enough, 
\begin{equation}
\label{eq:goal2}
    \max_{\xi\in [k,k+1]}F_\varepsilon(\xi)= F_\varepsilon(k+\frac12)=\frac{1}{\pi^4}, \mbox{ for all }k\in \Z,
\end{equation}
 where $$F_\varepsilon(\xi):=|\xi|^4 |\widehat{\tilde{f_2}}(\xi)|^2=\frac{1}{\pi^4}\sin^4(\pi \xi)+ \varepsilon ^2 g(\xi).$$

First, since $F_0$ is periodic and $F''_0(k+\frac12)<0$, by continuity and property (3), there exists a $\delta>0$ such that for $\varepsilon$ small enough  and all $k\in \Z$, \begin{equation}
\label{eq:deriv}
    F''_\varepsilon(\xi)<0\mbox{ for } |k+\frac12- \xi|<\delta.
\end{equation}

Using property (2), we see that $F_\varepsilon'(k+\frac12)=0$. Therefore, from \eqref{eq:deriv} we deduce that $$ \max_{\xi\in [k+\frac12- \delta,k+\frac12+\delta]}F_\varepsilon(\xi)= F_{\varepsilon}(k+\frac12)=\frac{1}{\pi^4}, \mbox{ for all }k\in \Z.$$

To complete the proof, we observe that  property (3) implies for $\varepsilon$ small enough, $$ \max_{\xi\in [k,k+1] \setminus [k+\frac12- \delta,k+\frac12+\delta]}F_\varepsilon(\xi)\leq \frac{1}{\pi^4}\cos^4( \pi \delta)+ \varepsilon^2 \norm{g}_\infty<\frac{1}{\pi^4},$$ that is, \eqref{eq:goal2} is true.

\end{proof}

Before giving the proof of item (3) in Theorem \ref{theorem:main}, we recall a well-known fact about the Fourier transform of radial functions on $\R^3$, see, e.g., \cite{grafakos}. 

\begin{lemma}
\label{lemma:g3}Let $F:\R^3 \to \mathbb{C}$ be a radial function, $F(x)=F_0(|x|)$. Then
    $$\widehat{F}(\xi)=\frac{2}{|\xi|}  \int_{0}^\infty F_0(x) {\sin({2 \pi |\xi|x}}) xdx. $$
    In particular,
    $$\widehat{\delta_{\lambda \mathbb{S}^{2}}}(\xi)= 2 \lambda  \frac{\sin{2 \pi \lambda  |\xi|}}{|\xi|}.$$
\end{lemma}

\begin{proof}[Proof of item (3) in Theorem \ref{theorem:main}]
Recall that here
 $\alpha=1,\beta=2$ and $d=3$.
    Set  $$f_3(x)=\mathbbm{1}_{[0,1]}(|x|) |x|^{-1}  \mbox{ and }H(x)= \psi_3(x)+ |x| -1,$$ where 
    $$\psi_3(x)=\frac{4}{\pi^3} \sum_{k=0}^\infty \frac{2}{(2k+1)^4} \widehat{\delta_{(k+\frac12)\mathbb{S}^2}}(x)=\frac{4}{\pi^3} \sum_{k=1}^\infty \frac{1-(-1)^k}{k^3} \frac{\sin{\pi k |x|}}{|x|}.$$

 Let us confirm that all the conditions of Lemma \ref{lemma:dual} hold.
First,  since

    $$\widehat{f_3}(\xi)=\frac{2}{|\xi|} \int_0 ^1 \sin{2 \pi |\xi| x} dx= \frac{1-\cos{2 \pi |\xi|}}{\pi |\xi|^2}$$ and $\norm{\widehat{f_3} |\xi|^2}_\infty=\frac2\pi$, we have that $\Lambda= \mathbb{N}+\frac12$. Thus, we see that both \eqref{eq:tekni} and property (iii) hold.

Second, to verify properties (i) and (ii) of Lemma \ref{lemma:dual}, we note  that, since
    $$\int_0 ^1 x(1-x) \sin{\pi k x}=\frac{2}{\pi^3k^3} (1-(-1)^k),$$  we have
    $$ x(1-x)=\phi(x),\qquad 0\leq x \leq 1,$$ where
$$\phi(x):=\frac{4}{\pi^3} \sum_{k=1}^\infty \frac{1-(-1)^k}{k^3} \sin{\pi k x},\quad x \in \R.$$

    Finally, using that $\phi(x+2)=\phi(x)$ and $\phi(x)=-\phi(-x)$, we see that      $$H(x)/|x|=\phi(x)+ |x|^2-|x|=\begin{cases}
        0, &0<|x|\leq 1\\
        2(|x|-1)^2, &1<|x|\le 2\\
       4|x|-6, &2<|x|<3.\\
        
    \end{cases}$$
Thus, $f_3 H\equiv 0$, and,
since $x^2-x$ is increasing for $x\geq \frac12$, we derive that $H\geq 0$.    

 It follows from a straightforward computation that
   $C^*_{1,2,3}(f_3)=
    \frac{4}{9 \pi ^2}.$

   The proof of the radial uniqueness is similar to that of item (2). By Lemma \ref{lemma:dual}, any radial extremizing $\mu\in \mathcal{M}^+(\R^d)$ normalized so that $\norm{\mu}_{TV}=\norm{f_3}_1$ and $\norm{\widehat{f_3}(\xi)|\xi|^2}_\infty=\norm{\widehat{\mu}(\xi)|\xi|^2}_\infty$ must be supported on the ball of radius 1 and, denoting the radial part of $\widehat{\mu}$ by $\widehat{\mu}_0$, 
    $$\widehat{\mu}_0(k+\frac12)=(\widehat{f_3})_0(k+\frac12), \; k \in \N \mbox{ and } \widehat{\mu}'_0(k+\frac12)=(\widehat{f_3})'_0(k+\frac12).$$ This means that the entire functions
    $$F(z)=\int_0^1x \sin({2 \pi z x} )(f_3)_0(x) dx \mbox{ and } G(z)=\int_0^1 x\sin({2 \pi z x}) d\mu_0(x) $$ satisfy $$F(k+\frac12)=G(k+\frac12) \mbox{ and } F'(k+\frac12)=G'(k+\frac12)$$ for $k \in \Z$. Thus, since $|F(z)|+|G(z)|\lesssim e^{2 \pi  |\Im{z}|}$, we conclude that $$g(z)=\frac{F(z)-G(z)}{\cos^2{ \pi z}}$$ is a bounded, entire function with $\liminf_{x \to \infty} |g(x)|=0$. This yields that $F=G$, that is, $\widehat{f_3}=\widehat{\mu}$.
\end{proof}

\begin{proof}[Proof of item (4) in Theorem \ref{theorem:main}]
Recall that by Lemma \ref{lemma:g3}
    $$\widehat{\mu_4}(\xi)=2 \frac{\sin{2 \pi  |\xi|}}{|\xi|}. $$
    Thus, for $0<\beta\leq 1$ $$\norm{\widehat{\mu_4}(\xi) |\xi|^\beta}_\infty=2 \frac{\sin{2 \pi  \lambda_\beta}}{\lambda_\beta^{1-\beta}}, $$
    where $\lambda_\beta$ is the unique zero of the derivative of $\widehat{\mu_4}(\xi) |\xi|^\beta$ in the interval $(0,\frac14]$. Equivalently,
    \begin{equation}
    \label{eq:deflb}
     2 \pi \lambda_\beta =(1-\beta)\tan{2 \pi \lambda_\beta }, \; \lambda_\beta \in (0,\frac14), 
    \end{equation} with the understanding that $\lambda_1=\frac14$.
    Set
    $$H(x)=A\widehat{\delta_{\lambda_\beta \mathbb{S}^2}}+ |x|^{\alpha}-B,$$ where
\begin{equation}
\label{eq:defAB}
    A=\frac{\alpha}{2 \lambda_\beta \beta \sin(2\pi \lambda_\beta)} \mbox{ and }B=\frac{\alpha}{\beta}+1.
\end{equation}

    We now verify the hypothesis of Lemma \ref{lemma:dual}. First, clearly \eqref{eq:tekni} and property (iii) hold. Second, for item (ii), note that
$$H(x)=2 \lambda_\beta A  \frac{\sin{2 \pi \lambda_\beta} |x|}{|x|}+ |x|^{\alpha}-B,$$ which readily shows that $H(x)=0$ for $|x|=1$. Finally, in Lemmas
\ref{lemma:inequality} and \ref{lemma:b1} it is shown that $
H(x)\geq 0 $ and $H(x)>0$ for $x\not \in \{0,1\}$,
 whence both item 
(i) of Lemma \ref{lemma:dual} and the uniqueness of $\mu_4$ in the class of radial measures follow.

Moreover, we find 

$$C_{\alpha,\beta,3}^*=C_{\alpha,1,3}(\mu_4)=\left(\frac{\sin{2 \pi \lambda_\beta}}{2\pi \lambda_\beta^{1- \beta}}\right)^\alpha.$$

\end{proof}

\vspace{0.2cm}\section{Duality and asymptotics of extremizers}
\label{sect:duality}
In this section we investigate the linear programming duality interpretation of inequality \eqref{uncer:measures}. As a result, we will obtain Theorem \ref{extremizer} and 
Corollary \ref{th:asymptotic}.

To this end, we define \begin{align*}
    E_1&=\left\{\mu \mbox{ radial probability measure such that } \norm{\widehat{\mu}(\xi)|\xi|^\beta}_\infty \leq 1 \right\}
\end{align*}
and observe that by Theorem \ref{theorem:notmain}  we have $$C^{*\frac{1}{\beta}}_{\alpha,\beta,d}=\min_{\mu\in E_1} I_{\mu}(|x|^\alpha).$$

The dual  of this minimization problem (see Chapter 5 in \cite{convex} for a detailed explanation) is given by $$\sup_{(\nu,\lambda)\in E_2} \big(\lambda- \norm{\nu}_{TV}\big),$$
where \begin{align*}
E_2&=\Big\{(\nu,\lambda) :  \nu \mbox{ radial finite signed measure }
\\&\qquad\qquad\qquad\quad
\mbox{and }
\lambda\in \R \mbox{ such that }  \widehat{|\xi|^\beta \nu} +|x|^\alpha -\lambda\geq 0\Big\}. 
\end{align*}
Here we note that the condition $$\widehat{|\xi|^\beta \nu} +|x|^\alpha -\lambda \geq 0$$ should be understood in the sense of distributions, that is,  $$0\leq \int g \left(\widehat{|\xi|^\beta \nu} +|x|^\alpha -\lambda\right), \qquad 0\leq g\in S(\R^d).
$$
The standard weak duality result (see
Step 1 in the proof of Lemma \ref{lemma:strongduality}) yields the inequality
$$C^{*\frac{1}{\beta}}_{\alpha,\beta,d}=\min_{\mu\in E_1} I_{\mu}(|x|^\alpha) \geq \sup_{(\nu,\lambda)\in E_2} \lambda- \norm{\nu}_{TV}.$$ As a consequence, if for some $\mu^*\in E_1$ and $(\nu^*,\lambda^*)\in E_2$ we have
$$ I_{\mu^*}(|x|^\alpha) = \lambda^*- \norm{\nu^*}_{TV},$$
then (cf. Lemma \ref{lemma:dual})
$$C^{*\frac{1}{\beta}}_{\alpha,\beta,d}= I_{\mu^*}(|x|^\alpha) = \lambda^*- \norm{\nu^*}_{TV}.$$

The goal of the next lemma is to show that strong duality holds and to obtain the corresponding slackness conditions. We make use of  
the Fenchel--Rockafellar duality result from
\cite{Optimaltransport}.


\begin{lemma}
\label{lemma:strongduality}
Let $\alpha,\beta>0$. Then 
\begin{equation}
 \label{eq:strongduality}
     C^{*\frac{1}{\beta}}_{\alpha,\beta,d}=\min_{\mu\in E_1} I_{\mu}(|x|^\alpha) = \max_{(\nu,\lambda)\in E_2} \lambda- \norm{\nu}_{TV}.
 \end{equation}

 Moreover, any extremizers $\mu^*$ and $(\nu^*,\lambda^*)$ satisfy the slackness conditions
   \begin{equation}
\label{eq:cond(i)}
 \int \widehat{\mu^*}(\xi) |\xi|^\beta d \nu^*(\xi)=  \norm{\nu^*}_{TV}
  \end{equation} 
and
        \begin{equation}
\label{eq:cond(ii)}\lim_n \int g_n \left(\widehat{|\xi|^\beta \nu^*} +|x|^\alpha -\lambda^*\right)= 0, \end{equation} 
where
        \begin{equation}
 \label{eq:gm}
     g_n(y)= \int_{\R^d} \phi_n(y-x)d \mu^* (x),
 \end{equation} for $\phi_n(x)=n^d \phi(nx)$ with $0\leq \phi\in S_{rad}(\R^d)$  and $\int_{\R^d} \phi=1$,
    where $S_{rad}(\R^d)$ is the space of radial Schwartz functions.
    
\end{lemma}

\begin{proof}

\textbf{{Step 1.}}
To obtain the weak duality inequality
 \begin{equation}
 \label{eq:sup}
    C^{*\frac{1}{\beta}}_{\alpha,\beta,d}\geq  \sup_{(\nu,\lambda)\in E_2} \lambda- \norm{\nu}_{TV},
 \end{equation} we first show that
 \begin{equation}
 \label{eq:inf}C_{\alpha,\beta,d}^{* \frac 1 \beta}=\inf_{\substack{\mu\in E_1\\
 d\mu=g dx,\, g\in S(\R^d)}} I_{\mu}(|x|^\alpha).
 \end{equation}
This follows by considering the convolutions \eqref{eq:gm}
 and noting that for any $\mu\in E_1$,
$$\norm{g_n}_1=\norm{\mu}_{TV}, \norm{\widehat{g_n}(\xi)|\xi|^\beta}_\infty \leq \norm{\widehat{\mu}(\xi)|\xi|^\beta}_\infty \mbox{ and } \lim_n \norm{g_n(x) |x|^\alpha}_1= I_{\mu}(|x|^\alpha).$$

Second, for any $(\nu,\lambda)\in E_2$ and $\mu\in E_1$ with $d \mu= g d x$ and $ g\in S(\R^d)$,

\begin{align*}
    0&\leq \int g \left(\widehat{|\xi|^\beta \nu} +|x|^\alpha -\lambda\right)\\
    &=-\lambda+ \norm{g(x)|x|^\alpha}_1 + \int \widehat{g}(\xi)|\xi|^\beta d \nu\\ &\leq -\lambda+ \norm{g(x)|x|^\alpha}_1 + \norm{\nu}_{TV}.
\end{align*} Hence, by using \eqref{eq:inf} we deduce that
 $$0\leq -\lambda+ C_{\alpha,\beta,d}^{* \frac 1\beta} + \norm{\nu}_{TV},$$ that is, \eqref{eq:sup} holds.

 \vspace{3mm}
\textbf{{Step 2.}}
In this part we prove \eqref{eq:strongduality}. For this,
it is sufficient  to obtain the strong duality inequality
 \begin{equation}
 \label{eq:max}
     C_{\alpha,\beta,d}^{* \frac 1 \beta}\leq  \max_{(\nu,\lambda)\in E_2} \lambda- \norm{\nu}_{TV}.
 \end{equation}
    Following Section 1.1.6 of \cite{Optimaltransport}, let $$E=\left\{g \mbox{ radial such that } |g(\xi)|(1+|\xi|^\beta)\in C_0 \right\}.$$ Then, the dual of $E$
    with respect to the duality pairing
    $$(g, \nu)= \int g d \nu,$$ 
    is given by $$E^*=\left\{ \nu \mbox{ radial such that}\int \frac{d|\nu|}{1+
    |\xi|^\beta}<\infty   \right\}.$$ 
    
    Define, interpreting the Fourier transform of $g\in E$ in the distributional sense,
$$     \Theta(g)=\begin{cases}
        0, &\norm{g(\xi)|\xi|^\beta}_\infty\leq 1\\
        +\infty, &\mbox{ otherwise}; \end{cases}
 $$   
 and
  $$  \Xi(g)=\begin{cases}
  \displaystyle \int |x|^\alpha \widehat{g}(x) dx, & \displaystyle\widehat{g}\geq 0, \int \widehat{g}=1\\
        +\infty, &\mbox{ otherwise}.
    \end{cases}
  $$
   Their Legendre-Fenchel transforms (see, e.g., \cite[p. 23]{Optimaltransport}) are given by
   \begin{align*}
       \Theta^*(\nu)&= \int |\xi|^{-\beta} d |\nu|,\\
    \Xi^*(\nu)&=\sup_{\substack{\widehat{g}\geq 0 \\\int \widehat{g}=1\\g \in E}} \int g d \nu - \int |x|^\alpha \widehat{g}(x) dx,
   \end{align*}
   respectively.
   Observe that both $\Theta$ and $\Xi$ are convex; and that $\Theta(g_0)$ and $\Xi(g_0)$ are finite for $g_0(x)= \phi(Nx)$ with $N$ large enough, where $0\leq \widehat{\phi} \in S_{rad}(\R^d)$ satisfies $\int\widehat{\phi}=1$. In addition, $\Theta$ is continuous at $g_0$ for $N$ large enough.
   
   Then, by the 
   Fenchel--Rockafellar duality theorem
(see  \cite[Theorem 1.9]{Optimaltransport}),
   \begin{equation}
       \label{eq:fenchel}
      C^{*\frac{1}{\beta}}_{\alpha,\beta,d} =\inf_{g\in E} \Big(\Theta(g)+ \Xi(g)\Big)= \max_{\nu \in E^*} \Big(-\Theta^*(-\nu)-\Xi^*(\nu)\Big).
   \end{equation} Let $\tilde{\nu}$ be an extremizer of the last maximum, then setting $$\lambda^*= -\Xi^*(\tilde{\nu}),$$
   we have that for any $0\leq \widehat{\phi} \in S_{rad}(\R^d)$ with $\int \widehat{\phi}=1$,
   $$\int \widehat{\phi} \left(- \widehat{ \tilde{\nu}}+ |x|^\alpha -\lambda^*\right)= -\int {\phi} d\tilde{\nu} +\int |x|^\alpha \widehat{\phi}(x) dx   - \lambda^* \geq 0. $$
   This means that $$-\widehat{ \tilde{\nu}}+ |x|^\alpha -\lambda^*\geq 0.$$ To finish the proof, recall that
   $$ C^{*\frac{1}{\beta}}_{\alpha,\beta,d} =-\int |\xi|^{-\beta} d |\tilde{\nu}| +\lambda^*.$$ Thus, setting $d\nu^*=-|\xi|^{-\beta}d\tilde{\nu}$, we have $$ C^{*\frac{1}{\beta}}_{\alpha,\beta,d} =-\norm{\nu^* }_{TV}+ \lambda^*$$ and $$\widehat{  |\xi|^\beta \nu^*}+ |x|^\alpha -\lambda^*\geq 0,$$
   that is, \eqref{eq:max} is valid.

 \vspace{3mm}
\textbf{Step 3.} We now obtain the slackness conditions.
Assume that $\mu^*$ and $(\nu^*,\lambda^*)$ satisfy  
\begin{equation}
    \label{eq:nomin}
     C^{*\frac{1}{\beta}}_{\alpha,\beta,d}= I_{\mu^*}(|x|^\alpha) = \lambda^*- \norm{\nu^*}_{TV}.
\end{equation}

To prove properties \eqref{eq:cond(i)} and 
\eqref{eq:cond(ii)}, consider the sequence $(g_n)_{n=1}^\infty $, where $g_n$ is defined in \eqref{eq:gm}. 
Then,
by the dominated convergence theorem, we obtain
  $$\lim_n  \int \widehat{g_n}(\xi) |\xi|^\beta d \nu^*= \int \widehat{\mu^*}(\xi) |\xi|^\beta d \nu^*.$$ Therefore,
\begin{align*}
    0&\leq \limsup_n \int g_n \left(\widehat{|\xi|^\beta \nu^*} +|x|^\alpha -\lambda^*\right) \\
    &\leq  -\lambda^* + \limsup_n \norm{g_n(x)|x|^\alpha}_1 + \limsup_n\int \widehat{g_n}(\xi) |\xi|^\beta d\nu^*\\
    &= -\lambda^* + C_{\alpha,\beta,d}^{* \frac 1 \beta} + \int \widehat{\mu^*}(\xi) |\xi|^\beta d \nu^*\\
    &\leq -\lambda^* +C_{\alpha,\beta,d}^{* \frac 1 \beta} + \norm{\nu^*}_{TV} =0.
\end{align*}

 Hence, we deduce that
  $$0= \limsup_n \int g_n \left(\widehat{|\xi|^\beta \nu^*} +|x|^\alpha -\lambda^*\right)$$
  and
  $$\int \widehat{\mu^*}(\xi) |\xi|^\beta d \nu^*=\norm{\nu^*}_{TV},$$ that is, \eqref{eq:cond(i)} and 
\eqref{eq:cond(ii)} hold.
\end{proof}

The proof of Theorem \ref{extremizer} is based on the following lemma:
\begin{lemma}
\label{th:structure}
Let $\mu\in \mathcal{M}^+(\R^d)$ be radial and satisfy $C_{\alpha,\beta,d}(\mu^*)=C^*_{\alpha,\beta,d}$. Then, one of the following must hold:
\begin{enumerate}[label=(\roman*)]
    \item for any $\Delta>0$ there exists $|\xi|>\Delta$ such that $$\norm{\widehat{\mu^*}(\xi)|\xi|^\beta}_\infty=|\widehat{\mu^*}(\xi)||\xi|^\beta;$$
    \item there exist $N,c_k,\lambda_k\geq 0$ such that \begin{equation}
    \label{mu* is finite2}
    \mu^*=\sum_{k=0}^N c_k\delta_{\lambda_k \mathbb{S}^{d-1}}.
\end{equation}

\end{enumerate}
\end{lemma}

\begin{proof} Assume that $\norm{\mu^*}_{TV}=\norm{\widehat{\mu^*}(\xi)|\xi|^\beta}_\infty=1.$
Let $$\Lambda:=\left\{|\xi|: |\widehat{\mu^*}(\xi)||\xi|^\beta =1\right\}.$$ We have to show that if  $\Lambda$ is bounded then \eqref{mu* is finite2} holds.

\textbf{Step 1.} By Lemma \ref{lemma:strongduality}, there exist $(\nu^*,\lambda^*)\in E_2$ such that
$$C_{\alpha,\beta,d}^{*\frac1 \beta}= I_{\mu^*}(|x|^\alpha)=\lambda^*- \norm{\nu^*}_{TV}.$$
We show that \begin{equation}
    \label{eq:Lambda}
    \operatorname{supp}\nu^*\subset\Lambda.
\end{equation}
   From \eqref{eq:cond(i)} in Lemma \ref{eq:strongduality}, we deduce that the pair $(\mu^*, \nu^*)$ satisfies $$\int d |\nu^*| - \int \widehat{\mu^*}(\xi)|\xi|^\beta d \nu^*=0.$$ Hence, if $|\widehat{\mu^*}(\xi_0)||\xi_0|^\beta<1=\norm{\widehat{\mu^*}(\xi)|\xi|^\beta}_\infty$ for some $\xi_0$, then by continuity there exist 
$\delta,\varepsilon>
0$ such that
$|\widehat{\mu^*}(\xi)||\xi|^\beta<1-\delta$ for any  $\xi\in B:=\{\xi:|\xi-\xi_0|<\varepsilon\}.$ Thus, $$\delta \int_B d|\nu^*|\leq  \int_B d|\nu^*| -  \int_{B} \widehat{\mu^*}(\xi)|\xi|^\beta d \nu^* \leq \int d |\nu^*| - \int \widehat{\mu^*}(\xi)|\xi|^\beta d \nu^*=0.$$ Therefore, $\xi_0 \not \in \operatorname{supp} \nu^* $ and \eqref{eq:Lambda} holds.

  \textbf{Step 2.}
   Thus, if $\Lambda$ is bounded, then $\nu^*$ is compactly supported. In particular, the distribution $\widehat{|\xi|^\beta \nu^*}$ is, in fact, an analytic function that is bounded on $\R^d$.
This implies that
  $$H(x):=\widehat{|\xi|^\beta \nu^*} +|x|^\alpha -\lambda^*$$ is a radial, non-negative, continuous function satisfying
   $\lim _{|x|\to \infty} H(x)=+\infty$.

\textbf{Step 3.}
We show that \begin{equation}
    \label{eq:Z}
    \operatorname{supp}{\mu^*}\subset\{H=0\}.
\end{equation} 

  Let $x_0$ satisfy $H(x_0)\neq0$. Then $H(x)>0$ for $|x-x_0|<\varepsilon$ and small $\varepsilon>0$. From \eqref{eq:cond(ii)} in Lemma \ref{lemma:strongduality} we deduce that, for any non-negative continuous function $\varphi$ supported on $|x-x_0|<\varepsilon$,
$$\int_{|x-x_0|<{\varepsilon}} H(x) \varphi d\mu^*=\lim_n \int    H \varphi g_n=0,$$ whence we see that $\mu^*(|x-x_0|<\varepsilon)=0$. 
 
\textbf{Step 4.} Finally, we show that $$\{H=0\}=\bigcup_{k=0}^N \lambda_k \mathbb{S}^{d-1},$$ for some $N,\lambda_k\geq 0.$ After establishing this, from the radiality of $\mu^*$ and the fact that
$$\operatorname{supp}{\mu^*}\subset \{H=0\}$$ we can conclude that \eqref{mu* is finite2} holds.

To obtain this result, observe that since $H$ is radial and $\lim_{|x|\to \infty} H(x)=+\infty$, it suffices to show that the zero set of the function $$H_0(t)=H(t,0,\dots,0), \; t \in \R$$  has no accumulation points.

From the fact that $\widehat{|\xi|^\beta \nu^*}$ is analytic, we deduce that $H_0$ is analytic on $\R\setminus \{0\}$. This shows that the only possible accumulation point of the zero set is $t=0$.

To rule out this possibility, observe that since $\left(\widehat{|\xi|^\beta \nu^*}\right)_0$ is not a polynomial, there exist $l,r\geq 1$ such that both $\left(\widehat{|\xi|^\beta \nu^*}\right)_0^{(l)}(0)$ and $\left(\widehat{|\xi|^\beta \nu^*}\right)_0^{(l+r)}(0)$ are nonzero.
Thus, for  $0<t \to 0$, 
$$H_0(t)=A + B t^l + C t^{l+r}  + O(t^{l+r+1})  + t^\alpha -\lambda^*, $$ with $B,C$ non-zero. 
 Hence,
 $$H_0(t)=A-\lambda^*+ D t^s + o(t^{s})\; \mbox{ for some } s \in \{l,\alpha,l+r\} \mbox{ and } D\neq 0,$$ which implies that the zeros of $H_0$ cannot accumulate at $t=0$.

\end{proof}
We are now in a position to prove Theorem \ref{extremizer} and Corollary \ref{th:asymptotic}.
Recall that  the Fourier transform of a  radial measure $\mu\in \mathcal{M}(\R^d)$ is given by (see the Appendix in \cite{grafakos})
\begin{equation}
    \label{eq:radialfourier}
    \widehat{\mu}(\xi)=\frac{2 \pi}{\omega_{d-1} |\xi|^{\frac{d}{2}-1}}\int J_{\frac d2 -1}(2 \pi |x \xi|) |x|^{1-\frac d2} d\mu(x),
\end{equation} where 
$J_k$
is the Bessel function of order $k$. Recall that  $J_k$ satisfies, for some $A_k,\theta_k \in \R_+$,
\begin{equation}
    \label{eq:bessel}
    J_k(r)=\begin{cases}
        A_k r^k + O(r^{k+1}), &\mbox{ as } r\to 0\\
        \sqrt{\frac{2}{\pi r}}  \cos({r-\theta_k}) +O(r^{-\frac32}), &\mbox{ as } r\to \infty.\\
    \end{cases}
\end{equation}

    \begin{proof}[Proof of Theorem \ref{extremizer} and Corollary \ref{th:asymptotic}]
We first show that if item (ii) in Lemma \ref{th:structure} holds, then
\begin{equation}
\label{eq:asymptgoal}
    0<\limsup_{|\xi|\to \infty} |\xi|^{\frac{d-1}{2}}\left|\widehat{\mu^*}(\xi)\right|<\infty.
\end{equation}
To see this, we use \eqref{eq:radialfourier} to deduce that
\begin{equation}
\label{eq:estbessel}
\begin{aligned}
    \widehat{\mu^*}(\xi)&= |\xi|^{-\frac{d}{2}+1}\sum_{k=0}^N  c'_k J_{\frac{d}{2}-1}(2 \pi \lambda_k |\xi|)\\
    &= |\xi|^{-\frac{d-1}{2}} \sum_{k=0}^N  c''_k \cos({\lambda_k |\xi|- \theta_k}) + O(|\xi|^{-\frac{d+1}{2}}),
    \end{aligned}
\end{equation}
for some $c'_k,c''_k > 0$.
Relation \eqref{eq:asymptgoal} follows from the fact that $$0<\limsup_{|\xi|\to\infty} \left|\sum_{k=0}^N  {c}''_k \cos({\lambda_k |\xi|- \theta_k}) \right|<\infty$$ (the first inequality can be seen, for instance, by using Ingham's inequality \cite[Theorem A]{ingham}).

Let now $\beta>\frac{d-1}{2}$. By using relation \eqref{eq:asymptgoal} we see that part (ii) in Lemma \eqref{th:structure} contradicts the condition $\norm{\widehat{\mu^*}(\xi) |\xi|^\beta}_ \infty<\infty$, which means that part (i) must hold. This proves Theorem \ref{extremizer} and Corollary \ref{th:asymptotic} for this regime of the parameter $\beta$. 

Let $0<\beta\leq \frac{d-1}{2}$. First, we prove that $\mu^*$ is compactly supported.  

 Let $(\nu^*,\lambda^*)\in E_2$ be from Lemma \ref{lemma:strongduality}. We are going to show that if $\phi \in S_{rad}(\R^d)$ is supported on $\{|x|\geq 1\}$, then
    \begin{equation}
    \label{eq:goalboundedness}
        \left |\int \widehat{|\xi|^\beta \nu^*} \phi \right|\leq C \norm{\phi}_1 \quad\mbox{with}\quad C<\infty.
    \end{equation} 
 From \eqref{eq:radialfourier} we have that

    $$ \left |\int \widehat{|\xi|^\beta \nu^*} \phi \right|\lesssim  \int |\phi(x)| |x|^{-\frac{d}{2}+1} \int   |J_{\frac d2 -1}(2 \pi |x \xi|)| |\xi|^{1-\frac d2 + \beta }  d|\nu^*|(\xi) dx. $$   
    Observe that from \eqref{eq:bessel} and $\beta\leq \frac{d-1}{2}$ we deduce that for $|x|\geq 1$
 \begin{align*}
      |x|^{-\frac{d}{2}+1}\int |J_{\frac d2 -1}(2 \pi |x \xi|)| |\xi|^{1-\frac d2 + \beta }  d|\nu^*|(\xi) &\lesssim \norm{\nu^*}_{TV}|x|^{-\beta}\\
      &\leq \norm{\nu^*}_{TV},  \end{align*}
whence \eqref{eq:goalboundedness} follows.

  Let $0\leq  \phi \in S_{rad}(\R^d)$ be supported on $\{|x|\geq R\}$ with $\phi(x)>0$ for $|x|> R$. From \eqref{eq:goalboundedness} we deduce that 
    \begin{equation}
    \label{eq:goal3}
        (R^\alpha - \lambda^* - C) \norm{\phi}_1\leq\int \phi  \left(\widehat{|\xi|^\beta \nu^*} +|x|^\alpha -\lambda^*\right).
    \end{equation} 
        Next, if there exists $x_0\in \operatorname{supp} \mu^*$ with $|x_0|>R$, then
    $$ \lim _n \norm{g_n \phi}_1= \lim _n \int g_n \phi =  \int \phi d \mu^*>0, $$ where $g_n$ is given by \eqref{eq:gm}. Hence, from \eqref{eq:goal3} applied to $g_n \phi$ and \eqref{eq:cond(ii)} in Lemma \ref{lemma:strongduality} we conclude that 
     $$ \lim _n \norm{g_n \phi}_1(R^\alpha - \lambda^* - C) \leq \lim _n \int g_n \phi  \left(\widehat{|\xi|^\beta \nu^*} +|x|^\alpha -\lambda^*\right)=0,$$ which implies that $R^\alpha - \lambda^* - C\le0$, that is, that $R$ is bounded.

  Finally, the rest of the  proofs of 
  Theorem \ref{extremizer} and Corollary \ref{th:asymptotic}  follow from Lemma \ref{th:structure} and \eqref{eq:asymptgoal}.
\end{proof}
   
\vspace{0.3cm}

\section{Further results}\label{sec5}
In this section, we present two further results 
derived from the techniques developed in the previous sections.

First, we find the extremizers of \eqref{uncer:measures} for $(\alpha,\beta,d)=(2,2,1)$ under the additional restriction that both $\mu$ and $\widehat{\mu}$ are non-negative, where $ \mu \in \mathcal{M}^+(\R)$ is even. This problem was considered in \cite{steinerberger2}.

\begin{proposition}\label{prop:positive fourier}Let 
$
f_2(x)=(1-|x|)\mathbbm{1}_{[-1,1]}(x)$.
Then $$\inf_{\substack{ 0\neq \mu \in \mathcal{M}^+(\R^d)\\ \widehat{\mu}\geq 0}} C_{2,2,1}(\mu)=C_{2,2,1}(f_2)=\frac{1}{(6 \pi)^2}.$$ Moreover, $f_2$ is the unique extremizer.
\end{proposition}
\begin{proof}
The proof is a minor modification of those of Theorem \ref{theorem:main}.
    Let $$H(x)=\psi_5(x)+3 x^2- 1,$$ where
    $$\psi_5(x)=\sum_{j=0}^\infty a_j \cos(\pi j x)$$
    with $$\; a_j= \int_{-1} ^1 (1-3x^2) \cos({\pi j x}) dx=\begin{cases}
        0, &j=0\\
      12 \frac{(-1)^{j+1}}{\pi^2 j^2}, &j>0.
        
    \end{cases}$$
 Note that since $\psi_5(x)=1-3x^2$ for $|x|\leq 1$ and $\psi_5$ is 2-periodic, we have that
    $$H(x)=\begin{cases}
        0,\; &0<x\leq 1\\
        12x-12,\; &1<x<2. 
    \end{cases}$$
    Therefore, using again that $\psi_5$ is $2$-periodic and observing that $3x^2-1$ is increasing for $x>0$, we deduce that $H(x)\geq 0$ for $x\geq 0$; and, by evenness, we conclude that $H\geq 0$.

    Finally, for any even $\mu \in \mathcal{M}^+(\R^d)$ such that $\widehat{\mu}\geq 0$ and normalized so that $\norm{\mu}_{TV}=\norm{f_2}_1$ and $\norm{\widehat{\mu} |\xi|^2}_\infty=\norm{\widehat{f_2} |\xi|^2}_\infty=\frac{1}{\pi^2}$, we have
    \begin{align*}
        0&\leq \int Hd \mu  = 3I_\mu(x^2)- \norm{\mu}_{TV} + \sum_{j=0}^\infty a_j \frac{\widehat{\mu}(j/2)+\widehat{\mu}(-j/2)}{2}\\
        &\leq 3I_\mu(x^2)- \norm{\mu}_{TV} + \sum_{j=0}^\infty a_{2j+1} \frac{\widehat{\mu }(j+1/2)+\widehat{\mu }(-j-1/2)}{2}\\
        &  \leq  3I_\mu(x^2)- \norm{\mu}_{TV} + \norm{\widehat{\mu}(\xi)|\xi|^\beta}_\infty\sum_{j=0}^\infty a_{2j+1} \big(j+\frac12\big)^{-2} \\
          &  =  3I_\mu(x^2)- \norm{f_2}_{1} + \norm{\widehat{f_2}(\xi)|\xi|^\beta}_\infty\sum_{j=0}^\infty a_{2j+1} \big(j+\frac12\big)^{-2},
    \end{align*} where the second inequality follows from $\widehat{\mu}\geq 0$. Since, for  $d\mu=f_2 dx$,
    every inequality becomes an equality
    (recall that $\widehat{f_2} (\xi)= \left(\frac{\sin(\pi \xi)}{\pi \xi}\right)^2$), we have that 
$$0= 3\norm{f_2(x)|x|^2}_{1}- \norm{f_2}_{1} + \norm{\widehat{f_2}(\xi)|\xi|^\beta}_\infty\sum_{j=0}^\infty a_{2j+1} \big(j+\frac12\big)^{-2}$$ whence we deduce that $\norm{f_2(x)|x|^2}_{1}\leq I_\mu(x^2),$ and $f_2$ is an extremizer.

To establish uniqueness, observe that for any extremizer $\mu$,
every inequali\-ty in the previous argument must, in fact, be an equality.
Thus, $$\operatorname{supp} \mu\subset[-1,1],\; \widehat{\mu}(j)=\delta_0(j) \norm{f_2}_1, j \in \Z \mbox{ and } \widehat{\mu}(j+\frac12)=\frac{1}{\pi^2(j+\frac12)^2}, j \in \Z,  $$ which implies that $d\mu=f_2 dx$. This can be seen by considering the entire function
$$g(z)=\frac{\widehat{\mu}(z)-\widehat{f_2}(z)}{\sin{2 \pi z}},$$ which is bounded and satisfies $\liminf_{x\to \infty} |g(x)|=0$.
\end{proof}
Second, in some cases we can strengthen Theorem \ref{extremizer}  to show that the Fourier transform of any extremizer $\mu^*$ must oscillate infinitely often between $\displaystyle-\frac{\|{\widehat{\mu^*}(\xi) |\xi|^\beta}\|_\infty}{|\xi|^\beta} $ and  $\displaystyle\frac{\|{\widehat{\mu^*}(\xi) |\xi|^\beta}\|_\infty}{|\xi|^\beta} $. More precisely,
\begin{proposition}
    Let $\alpha \in 2 \mathbb{N}$ and let $\mu^* \in \mathcal{M}^+(\R^d)$ be radial and satisfy $C^{*}_{\alpha,\beta,d}=C_{\alpha,\beta,d}(\mu^*)$. Assume that there exist $\varepsilon,\delta>0$ such that $d\mu^*\geq \varepsilon dx$ on $\{ x:|x|<\delta\}$.  Define, for $\varepsilon \in \{+,-\}$,
     $$\Lambda_\varepsilon:=\Big\{|\xi|: \widehat{\mu^*}(\xi)|\xi|^\beta =\varepsilon\norm{\widehat{\mu^*}(\xi)|\xi|^\beta}_\infty\Big\}.$$ Then both $\Lambda_+$ and $\Lambda_-$ are unbounded.
\end{proposition}
\begin{proof}  
    For the sake of contradiction, assume that $\Lambda_-$ is bounded. Let $(\nu^*, \lambda^*)$ be those given in Lemma \ref{lemma:strongduality}. Write $\nu^*=\nu_1-\nu_2$ with $\nu_i \in \mathcal{M}^+(\R^d)$. Then, the proof of Step 1 in Lemma \ref{th:structure} shows that
    $$\operatorname{supp} \nu_1 \subset \Lambda_+ \mbox{ and } \operatorname{supp} \nu_2 \subset \Lambda_-,$$ and, in particular, we have that $\operatorname{supp} \nu_2$ is bounded.
    
    Next, we observe that from \eqref{eq:cond(ii)} in Lemma \ref{lemma:strongduality} and the condition $d\mu^*\geq \varepsilon dx$,  for any $0\leq \varphi\in S(\R^d)$ supported on $\{x:|x|<\delta\}$
    \begin{equation}
    \label{eq:thing}
        \int \varphi \left(\widehat{|\xi|^\beta \nu^*} +|x|^\alpha -\lambda^*\right)=0. 
    \end{equation}
    Indeed, the hypothesis  $d\mu^*\geq \varepsilon dx$ on $\{ x:|x|<\delta\}$ implies that, for $n$ large enough, the $g_n$ of \eqref{eq:cond(ii)} in Lemma \ref{lemma:strongduality} satisfy $g_n(x)\gtrsim 1 $ for $|x|\leq \delta$. Thus, there exists $C>0$ such that $g_n-C \varphi\geq 0$, and, by positivity, 
    
 $$0\leq \int \left(g_n-C\varphi \right)\left(\widehat{|\xi|^\beta \nu^*} +|x|^\alpha -\lambda^*\right), $$ whence \eqref{eq:thing} follows by letting $n\to \infty$ and using item (ii) in Lemma \ref{th:structure}.

    Thus, defining $\varphi_j(x)=j^d \varphi(jx),$ and additionally assuming that $\widehat{\varphi}\geq 0$ and $\widehat{\varphi}(0)=1$, we have
    \begin{align*}
       0&\leq \int \widehat{\varphi}(\xi/j) |\xi|^\beta d \nu_1 \\
       &= \int \widehat{\varphi}(\xi/j) |\xi|^\beta d \nu_2 + \int (\lambda^*-|x|^\alpha) \varphi_j \lesssim 1 ,
    \end{align*} where the equality follows from \eqref{eq:thing}; and the second inequality, from the fact that both $\nu_2$ and $\varphi$ are compactly supported.
    Hence, by Fatou's Lemma we deduce that $\int |\xi|^\beta d\nu_1<\infty$. This means that the continuous function $\widehat{|\xi|^\beta \nu_1}$ satisfies
    $$\widehat{|\xi|^\beta \nu_1}=\widehat{|\xi|^\beta \nu_2} -|x|^{\alpha} + \lambda^* \mbox{ for } |x|<\delta.$$ 
    Hence, using that $\alpha\in 2\N$, we see that  $\widehat{|\xi|^\beta \nu_1}$ is analytic on a neighborhood of the origin. Since $|\xi|^\beta \nu_1 \in \mathcal{M}^+(\R^d)$, by the propagation of regularity phenomenon (see Theorem 2.3 in \cite{propagation}), we have that $\widehat{|\xi|^\beta \nu_1}$ is analytic on $\R^d$. Therefore,
    $$\widehat{|\xi|^\beta \nu^*}= -|x|^{\alpha} + \lambda^* \mbox{ for all } x,$$ which is impossible because both $|\xi|^\beta \nu_1$ and $|\xi|^\beta \nu_2$ are finite measures, and, in particular, $\widehat{|\xi|^\beta \nu^*}$ is a bounded function.
\end{proof}
We note that the conditions $\alpha \in 2\N$ and $\mu^* \geq \varepsilon$ on $\{ x:|x|<\delta\}$ cannot be omitted. This follows  from 
 item (4) of Theorem \ref{theorem:main}.

\section{Appendix: trigonometric inequalities}
\label{section:auxi}
Define
$$H(x)=2 \lambda_\beta A  \frac{\sin{2 \pi \lambda_\beta} |x|}{|x|}+ |x|^{\alpha}-B,$$
where (see \eqref{eq:defAB})
$$    A=\frac{\alpha}{2 \lambda_\beta \beta \sin(2\pi \lambda_\beta)} \mbox{ and }B=\frac{\alpha}{\beta}+1.
$$

\begin{lemma}
\label{lemma:inequality}
 Let $\beta=1$ and $2 \geq \alpha \geq \alpha_0:=\frac{1}{\frac{\pi}{2}-1}$. Then, for any $x\geq 0$, we have $H(x)\geq 0,$ and $H(x)>0$ for $|x|\not \in \{0,1\}$.
    
\end{lemma}
    
\begin{proof}

Noting that  $\lambda_1=\frac 14$, 
 we define
 $$G(x):=xH_0(x)=2 \lambda_\beta A \sin{(2 \pi \lambda_\beta x)} +x^{\alpha+1}- Bx=\alpha \sin{\frac{\pi}{2} x}+ x^{\alpha+1}-(\alpha+1)x.$$ 
    We observe that 
    \begin{align*}
        G^{(1)}(x)=&\alpha \frac{\pi}{2} \cos{\frac{\pi}{2}x}+ (\alpha+1)x^{\alpha}- (\alpha+1)\\
        G^{(2)}(x)=&-\alpha \frac{\pi^2}{4} \sin{\frac{\pi}{2}x}+ \alpha (\alpha+1)x^{\alpha-1}.
    \end{align*}
Here note that $G^{(2)}(x)>0$ for $x\geq 1$, because for $\alpha\geq \alpha_0$ we have $\alpha (\alpha+1) >\alpha \frac{\pi^2}{4}$. Moreover, since $G(1)=G^{(1)}(1)=0$, we deduce that $G(x)> 0$ for $x> 1$.

 We prove now that $G(x)> 0$  for $x\in (0,1)$. Observe that for all $x\in [0,1]$,
$G^{(5)}(x) \geq 0.$
By using Taylor expansions as $x\to 0$, we derive 
\begin{equation}
\label{eq:taylor}
\begin{aligned}
    G(x)-x^{\alpha+1}=& \left(\alpha \frac{\pi}{2} - (\alpha+1)\right)x - \alpha \left(\frac{\pi}{2}\right)^3 \frac{x^3}{3!}  +\alpha\left(\frac{\pi}{2}\right)^5 \frac{x^5}{5!} +o(x^5), \\
    G(1-x)=&\alpha \cos {\frac{\pi}{2}x} + (1-x)^{\alpha+1}-(\alpha+1)(1-x) \\ =&\left( (\alpha+1)\alpha - \alpha \left( \frac{\pi}{2}\right)^2\right) \frac{x^2}{2} - (\alpha+1)\alpha (\alpha-1)\frac{x^3}{3!}\\
    &+\left( \alpha \left(\frac{\pi}{2}\right)^4+ (\alpha+1)\alpha(\alpha-1)(\alpha-2)\right)\frac{x^4}{4!}\\
    &- (\alpha+1)\alpha(\alpha-1)(\alpha-2)(\alpha-3)\frac{x^5}{5!}
    +o(x^5). 
    \end{aligned}
\end{equation}Then noting that
$\left(G(1-\cdot)\right)^{(k)}(x)=(-1)^k G^{(k)}(1-x)$,
we see that the following hold for $x> 0$ small enough:
\begin{enumerate}
   
    \item $G^{(4)}(x)< 0$ and $G^{(4)}(1-x)> 0$;
    \item $G^{(3)}(x)> 0$ and $G^{(3)}(1-x)> 0$;
    \item $G^{(2)}(x)> 0$ and $G^{(2)}(1-x)> 0$;
    \item $G^{(1)}(x)> 0$ and $G^{(1)}(1-x)< 0$ (since $\alpha_0\leq \alpha$);
    \item $G(x)> 0$ and $G(1-x)>0$.
\end{enumerate}
The fact that $G\geq 0$ now follows from
\begin{enumerate}[label=(\alph*)]
    \item Since $G^{(5)}(x)\geq 0$ for $x\in [0,1]$, from (1) we deduce that $G^{(4)}(x)\leq 0$ for $x\in [0,x_0]$ and $G^{(4)}(x)\geq 0$ for $x\in [x_0,1]$.
\item From (2) and (a), we see that either $G^{(3)}(x) \geq 0$ for all $x\in[0,1]$ or $G^{(3)}(x)\geq 0$ for $x\in [0,x_1] \cup [x_2,1]$ and $G^{(3)}(x)\leq  0$ for $x\in [x_1,x_2]$.
\item If $G^{(3)}(x)\geq 0$ for $x\in [0,1]$, from (3) we deduce that $G^{(2)}(x)\geq 0$ for $x\in [0,1]$, but this contradicts (4).

\item From (b), (c) and (3), we conclude  that either $G^{(2)}(x)\geq 0$ for all $x\in[0,1]$ (which is again not possible by (c)) or $G^{(2)}(x)\geq 0$ for $x\in [0,x_3] \cup [x_4,1]$ and $G^{(2)}(x)\leq  0$ for $x\in [x_3,x_4]$.

\item From (d) and (4), we deduce that $G^{(1)}(x)\geq 0$ for $x\in [0,x_5]$ and $G^{(1)}(x)\leq 0$ for $x\in [x_5,1]$.

\item From (e) and (5), we deduce that $G(x)> 0$ for $x\in (0,1)$.
\end{enumerate}
\end{proof}
For the proof of the case $0<\beta\leq 1$ we need an auxiliary \begin{lemma}
\label{lemma:inequalitiestan}
    $$1-\frac{x^2}{3}\geq \frac{x}{\tan{x}}\geq 2 \sqrt{1-\frac{x^2}{3}}-1,\;  x \in [0,\frac \pi 2].$$
\end{lemma}
\begin{proof}
    We use the power expansion of $x/\tan{x}$ at zero, which is valid for $|x| < {\pi}$,
    $$\frac{x}{\tan{x}}= \sum_{n= 0}^\infty \frac { (-1)^n 2^{2 n} B_{2 n} \, x^{2 n } } { (2 n)!}\leq 1-\frac{x^2}{3},$$ where we took into account  the facts that the Bernoulli numbers satisfy $B_{2n}(-1)^n\leq 0$ for $n\geq 1$, $B_0=1$, and $B_2=\frac{1}{6}.$
    The second inequality follows by using the expansion
    $$2 \sqrt{1- \frac{x^2}{3}} -1 = 1+  2\sum_{n=1}^\infty \binom{\frac12}{n} \left(\frac{-x^2}{3}\right)^n$$ and noting that 
    $$2 \binom{\frac12}{n} \left(\frac{-1}{3}\right)^n \leq \frac { (-1)^n 2^{2 n} B_{2 n} \,  } { (2 n)!},\; n \geq 1,$$ or, equivalently,
    $$2 \left|\binom{\frac12}{n}\right| \left(\frac{1}{3}\right)^n \geq \frac {  2^{2 n} \left|B_{2 n}\right| \,  } { (2 n)!},\; n \geq 1.$$

    For $n=1,2$ the latter can be verified directly.
    To prove it for $n\geq 3$, we use Theorem 1.1 in \cite{qi}, which gives

    $$\left|\frac{B_{2n+2}}{B_{2n}}\right|<\frac{1}{4 \pi^2 } (2n+1)(2n+2).$$
    Thus,
    $$\frac{\frac {  2^{2 n+2} \left|B_{2 n+2}\right| \,  } { (2 n+2)!}}{\frac {  2^{2 n} \left|B_{2 n}\right| \,  } { (2 n)!}} < \frac{1}{\pi^2}$$ which is less than
    $$\frac{2 \left|\binom{\frac12}{n+1}\right| \left(\frac{1}{3}\right)^{n+1} }{2 \left|\binom{\frac12}{n}\right| \left(\frac{1}{3}\right)^n }=\frac{1}{3} \frac{n-\frac12}{n+1}\geq \frac{1}{6} $$ for $n\geq 2$.
\end{proof}

 \begin{lemma}
\label{lemma:b1}
Let $0<\beta\leq 1$ and $\alpha\geq 2$. Then, for any $x\geq 0$, we have $H(x)\geq 0,$ and $H(x)>0$ for $|x| \neq 1$.
    \end{lemma}
\begin{proof}
As before, we define
$$G(x)=xH_0(x)=2 \lambda_\beta A \sin{(2 \pi \lambda_\beta x)} +x^{\alpha+1}- Bx$$ and, by using 
\eqref{eq:deflb}
and \eqref{eq:defAB}, we observe that $G(1)=G'(1)=0$.

First, we need the estimate
\begin{equation}
    \label{eq:boundslambdab}
    3 \beta > \left(2 \pi  \lambda_\beta\right)^2> 3 \beta-\frac{3}{4} \beta^2, \,0<\beta\leq 1,
\end{equation} which we now proceed to prove.
The upper bound follows by using \eqref{eq:deflb} and the upper bound in Lemma \ref{lemma:inequalitiestan}: $$1-\frac{(2 \pi \lambda_\beta )^2}{3} >\frac{2 \pi \lambda_\beta }{\tan{2 \pi \lambda_\beta }} =1-\beta.$$
To obtain the lower bound in \eqref{eq:boundslambdab}, we instead use the lower bound in Lemma \ref{lemma:inequalitiestan}:
$$1-\beta=\frac{2 \pi \lambda_\beta}{\tan{ 2 \pi \lambda_\beta}} > 2 \sqrt{1- \frac{(2 \pi \lambda_\beta)^2}{3}} -1,\quad 0<\beta\le 1.$$ 

Second, making use of \eqref{eq:boundslambdab}, we obtain $G^{(1)}(0)>0$ and $G^{(2)}(1)>0$. Indeed,
 $$G^{(1)}(0)= 4 \pi   A \lambda_\beta^2 -B= \frac{\alpha}{\beta}\left(\frac{2 \pi \lambda_\beta}{\sin{ 2 \pi \lambda_\beta}} - 1-\frac{\beta}{\alpha}\right)>0;$$
 because, taking into account \eqref{eq:boundslambdab},
 $$\left(\frac{2 \pi \lambda_\beta}{\sin{ 2 \pi \lambda_\beta}}\right)^2= \left(\frac{2 \pi \lambda_\beta}{\tan{ 2 \pi \lambda_\beta}}\right)^2 + \left(2 \pi \lambda_\beta\right)^2= (1-\beta)^2+\left({2 \pi \lambda_\beta}\right)^2> (1+ \frac{\beta}{2})^2.$$ 
 Again, by \eqref{eq:boundslambdab}, 
$$G^{(2)}(1)=-\frac{\alpha}{\beta} (2 \pi \lambda_\beta)^2+(\alpha+1)\alpha >0.$$

 Third, we prove that $G(x)>0$ for $x>1$. To see this, note that for $x\geq 1$
    $     G^{(2)}(x)\geq 0$. Indeed, in light of the inequality 
    $\frac{ \sin{ 2 \pi \lambda_\beta}}{2 \pi \lambda_\beta} \geq \frac{ \sin{ 2 \pi \lambda_\beta x }}{2 \pi \lambda_\beta x},\, x\geq 1, $
    we have
\begin{align*}
    G^{(2)}(x)&= -  \left(2 \pi \lambda_{\beta}\right)^2 2A \lambda_\beta   \sin ({2 \pi \lambda_\beta x}) +(\alpha+1)\alpha x^{\alpha-1}\\
    &\geq  -  \left(2 \pi \lambda_{\beta}\right)^2 2A \lambda_\beta x\sin ({2 \pi \lambda_\beta }) +(\alpha+1)\alpha x^{\alpha-1}\\
    &=-\frac{\alpha}{\beta} (2 \pi \lambda_\beta)^2x+(\alpha+1)\alpha x ^{\alpha-1}>0.
\end{align*}
    Thus, since $G(1)=G'(1)=0$, we deduce that $G(x)\geq 0$ for $x\geq 1$ and $G(x)=0$ only for $x=1$.
    
Fourth, we prove that $G(x)>0$ for $0< x <1$.
\begin{enumerate}[label=(\alph*)]
\item From $$G^{(3)}(x)=-(2 \pi \lambda_\beta)^3 2 \lambda_\beta A \cos ({2 \pi \lambda_\beta x}) +(\alpha-1)\alpha(\alpha+1)x^{\alpha-2}$$ and $0<\lambda_\beta\leq \frac14,$ we see that $G^{(3)}$ is non-decreasing on $[0,1]$.
\item   Since $G^{(3)}(0)\leq 0$,  we have that $G^{(3)}(x)\leq 0$ for $x\in [0, x_0]$ and $G^{(3)}(x)\geq 0$ for $x\in [x_0,1]$. The only other possibility is $G^{(3)}\leq 0$ on $[0,1]$, which is impossible because $G^{(2)}(0)=0$ and $G^{(2)}(1)>0$.
   \item  Since $G^{(2)}(0)=0$ and $G^{(2)}(1)>0$,
    we deduce that $G^{(2)}(x)\leq 0$ for $x\in [0, x_1]$ and $G^{(2)}(x)\geq 0$ for $x\in [x_1,1]$. 
 \item From $G^{(1)}(0)>0$ and $G^{(1)}(1)=0$ we see that
    $G^{(1)}(x)> 0$ for $x\in [0, x_2]$ and $G^{(1)}(x)\leq 0$ for $x\in [x_2,1]$.
    \item From $G(0)=G(1)=0$, we finally conclude that $G(x)> 0$ for $0< x < 1$.
\end{enumerate}
\end{proof}
\textbf{Acknowledgements.} We thank Egor Kosov for his helpful suggestions and for bringing reference \cite{Optimaltransport} to our attention.

\bibliographystyle{abbrv}
\bibliography{references.bib}

\end{document}